\input amstex
\documentstyle{amsppt}

\nologo
\magnification=1200
\overfullrule=0pt
\hsize = 5.5 true in\hoffset=.5in
\vsize = 7.5 true in
\loadeusm

\let\qed=\relax


\def\fn{{\frak n}}

\def\fnl{{\fn_L}}
\def\dfnl{{\dot\fn_L}}
\def\fnlq{{\fn_L^Q}}

\def\fnn{{\fn_N}}

\def\fns{{\fn_\fs}}
\def\fnz{{\fn_Z}}

\def\fm{{\frak m}}

\def\fS*{{\frak S_*}}
\def\fs{{\frak s}}

\def\ds{{\dot s}}
\def\dst{{\ds_T}}

\def\dstp{{\ds_{T'}}}
\def\dt{{\dot t}}

\def\dx{{\dot x}}
\def\dy{{\dot y}}

\def\fW{{\frak W}}

\def\sA{{\eusm A}}
\def\cA{{\Cal A}}
\def\sB{{\eusm B}}
\def\cB{{\Cal B}}
\def\sC{{\eusm C}}
\def\cC{{\Cal C}}
\def\sD{{\eusm D}}

\def\cE{{\Cal E}}

\def\cF{{\Cal  F}}

\def\cG{{\Cal G}}

\def\cH{{\Cal H}}
\def\wcH{{\widetilde{\Cal H}}}

\def\sM{{\eusm M}}
\def\cM{{\Cal M}}
\def\tM{{\widetilde {\sM}}}
\def\sN{{\eusm N}}
\def\cN{{\Cal N}}

\def\sU{{\eusm U}}

\def\sR{{\eusm R}}

\def\R{{\Bbb R}}
\def\T{{\Bbb T}}

\def\Z{{\Bbb Z}}

\def\a{{\alpha}}

\def\da{{\dot\a}}
\def\b{{\beta}}

\def\f{{\varphi}}
\def\p{{\psi}}
\def\z{{\zeta}}
\def\Ad{{\text{\rm Ad}}}

\def\tilh{{\tilde h}}

\def\tp{{\tilde p}}
\def\bp{{\bar p}}
\def\tq{{\tilde  q}}
\def\bq{{\bar q}}
\def\tr{{\tilde r}}
\def\br{{\bar r}}

\def\txti{{\text{\rm i}}}
\def\txk{{\text{\rm k}}}
\def\txm{{\text{\rm m}}}

\def\txmthth1{{\txm_{\thth1}}}
\def\txn{{\text{\rm n}}}

\def\txnn{{\txn_N}}

\def\txr{{\text{\rm r}}}
\def\txs{{\text{\rm s}}}

\def\wtQ{{\widetilde Q}}

\def\Qm{{Q_\txm}}

\def\wtG{{\widetilde G}}
\def\wtH{{\widetilde H}}
\def\wcH{{\widetilde {\Cal H}}}

\def\th{{\theta}}
\def\tht{{\theta_T}}
\def\tht'{{\th_{T'}}}
\def\la{{\lambda}}

\def\AFD{{\text{\rm AFD}}}
\def\Aut{{\text{\rm Aut}}}
\def\Autf'{{\Aut_\f'}}
\def\Autp'{{\Aut_\p'}}
\def\Cntp'{{\Cnt_\p'}}
\def\Cnt{{\text{\rm Cnt}}}
\def\Ob{{\text{\rm Ob}}}
\def\Obm{{\Ob_{\text{\rm m}}}}
\def\cnt{{\Cnt}}
\def\cntr{{\Cnt_{\text{\rm r}}}}
\def\Hom{{\text{\rm Hom}}}

\def\Iso{{\text{\rm Iso}}}

\def\Int{{\text{\rm Int}}}
\def\Intp'{{\Int_\p'}}
\def\Intb{{\overline{\text{\rm Int}}}}
\def\Ind{{\text{\rm Ind}}}
\def\Ker{{\text{\rm Ker}}}

\def\mod{{\text{\rm mod}}}

\def\Out{{\text{\rm Out}}}
\def\out{{\text{\rm out}}}

\def\Outt{{\Out_{\tau, \th}}}

\def\Autf{{\Aut_\f}}
\def\Autfm{{\Aut_\f(\sM)}}

\def\Autf'{{\Aut_\f'}}

\def\two{{\rm I\!I}}
\def\twoone{{\rm I\!I$_1$}}
\def\threee{{\text{\rm I\!I\!I}}}
\def\threeone{{\rm I\!I\!I$_1$}}
\def\threel{{\rm I\!I\!I$_{\lambda}$}}
\def\three0{{\rm I\!I\!I$_0$}}

\def\twoinf{{\rm I\!I$_{\infty}$}}

\def\tsU{{\widetilde {\sU}}}

\def\tAd{{\widetilde {\text{\rm Ad}}}}
\def\ta{{\widetilde {\a}}}

\def\d{{\delta}}

\def\dnu{{\dot \nu}}

\def\tla{{\tilde \la}}

\def\tmu{{\tilde \mu}}

\def\La{{\Lambda}}
\def\g{{\gamma}}

\def\bg{{\bar g}}
\def\bh{{\bar h}}

\def\tilh{{\tilde h}}

\def\tchi{{\widetilde {\chi}}}

\def\part{{\partial}}

\def\partth{{\part_\th}}

\def\log{{\text{\rm log}}}

\def\sig{{\sigma}}

\def\dsig{{\dot\sig}}
\def\sigf{{\sigma^{\f}}}

\def\sigfs{{\sigma_s^{\f}}}

\def\sigp{{\sigma^{\psi}}}

\def\sigps{{\sigma_s^{\psi}}}

\def\sigpdsg{{\sig_{\ds(g)}^\p}}
\def\sigpdsh{{\sig_{\ds(h)}^\p}}
\def\sigpdsgh{{\sig_{\ds(gh)}^\p}}
\def\sigpdsgph{{\sig_{\ds(g)+\ds(h)}^\p}}

\def\sigrs{{\sig_s^\rho}}

\def\dsigs{{\dsig_s}}

\def\wt{{semi-finite normal weight}}
\def\fwt{{faithful \wt}}
\def\botimes{{\bar \otimes}}
\def\id{{\text{\rm id}}}

\def\r0{{\sR_0}}

\def\r01{{\sR_{0,1}}}

\def\Map{{\text{\rm Map}}}

\def\botimes{{\overline \otimes}}
\def\wt{{semi-finite normal weight}}

\def\fwt{{faithful \wt}}

\def\vna{{von Neumann algebra}}

\def\tB{{\text{\rm B}}}

\def\tC{{\text{\rm C}}}

\def\tH{{\text{\rm H}}}

\def\tZ{{\text{\rm Z}}}

\def\tbsthr{{\tB_\txs^3}}
\def\thsthr{{\tH_\txs^3}}
\def\tzsthr{{\tZ_\txs^3}}

\def\dprime{{^{\prime\prime}}}

\def\inv{{^{-1}}}

\document

\def\X{{\text{\rm Xext}}}

\def\Res{{\text{\rm Res}}}

\def\inf{{\text{\rm inf}}}

\def\dfs{{\dot\fs}}

\def\sh{{\fs\!_{\scriptscriptstyle H}}}

\def\sp{{\fs_\pi}}

\def\sj{{\fs_j}}

\def\pig{{\pi\!_{\scriptscriptstyle G}}}

\def\piz{{\pi\!_{\scriptscriptstyle\tZ}}}
\def\dpi{{\dot \pi}}
\def\pim{{\pi_{\txm}}}

\def\pims{{\pi_{\txm}^*}}

\def\sz{{\fs_{\!\scriptscriptstyle\tZ}}}

\def\etaa{{\eta_\a}}

\def\etaq{{\eta_Q}}

\def\etain{{\eta_{\text{in}}}}
\def\etat{{{\eta_{T}}}}

\def\etat'{{\eta_{T'}}}

\def\tpi{{\tilde \pi}}

\def\_#1{{_{#1}}}
\def\^#1{{^{#1}}}

\def\ainv#1{{\a_{#1}^{-1}}}

\def\as#1#2{{\a_{\scriptscriptstyle(#1)}^{#2}}}
\def\bs#1#2{{\b_{\scriptscriptstyle(#1)}^{#2}}}

\def\(#1){{^{({#1})}}}

\def\scirc{{\lower-.3ex\hbox{{$\scriptscriptstyle\circ$}}}}
\def\bracet#1#2{{\{#1\}_{#2}}}
\def\bracettt#1{{\left\{#1\right\}_T}}
\def\bracett'#1{{\left\{#1\right\}_{T'}}}
\def\bracet'm#1{{\bracett'{\txm\left({#1}\right)}}}

\def\brackett#1{{\left[{#1}\right]_{T}}}
\def\brackett'#1{{\left[{#1}\right]_{T'}}}
\def\rtz{{\R/T\Z}}
\def\rt'z{{\R/T'\Z}}
\def\coron{{\!:\ }}
\def\cdg{{\text{countable discrete group}}}
\def\dg{{discrete group}}
\def\cdag{{\text{countable discrete amenable group}}}

\def\Sub{{\text{\rm Sub}}}
\def\hjr{{{\text{\rm HJR}-exact sequence}}}
\def\mhjr{{\text modified \hjr}}

\def\tcatw{{\tC_\a^2}}
\def\tzatw{{\tZ_\a^2}}

\def\tbmsout{{\tB_{\txm, \fs}^\out}}

\def\tzath{{\tZ_\a^3}}
\def\tzasth{{\tZ_{\a, \text{\rm s}}^3}}

\def\thasth{{\tH_{\a, \text{\rm s}}^3}}

\def\tbth1{{\tB_\th^1}}

\def\tbasout{{\tB_{\a, \fs}^\out}}
\def\tbmsout{{\tB_{\txm, \fs}^\out}}

\def\tctw{{\tC^2}}
\def\tzth1{{\tZ_\th^1}}
\def\thth1{{\tH_\th^1}}
\def\ththr{{\tH^3}}
\def\thsthr{{\tH_\txs^3}}
\def\ththra{{\tH_\a^3}}

\def\thaout{{\tH_\a^\out}}

\def\thsout{{\tH_{\fs}^\out}}
\def\thasout{{\tH_{\a, \fs}^\out}}

\def\thmsout{{\tH_{\txm, \fs}^\out}}

\def\tzthr{{\tZ^3}}

\def\tbthr{{\tB^3}}

\def\thtw{{\tH^2}}

\def\tztw{{\tZ^2}}
\def\tzatw{{\tZ_\a^2}}

\def\tzaout{{\tZ_\a^\out}}

\def\tzasout{{\tZ_{\a, \fs}^\out}}

\def\Grp{{\text{\rm Grp}}}
\def\chim{{\chi_{\text{\rm m}}}}

\def\pr{{\text{\rm pr}}}

\def\pijj'{{\pi_{J, J'}}}
\def\pijj's{{\pi_{J, J'}^*}}

\def\tmu{{\tilde \mu}}

\def\linfrt'z{{L^\infty(\R/T'\Z)}}

\def\rtz{{\R/T\Z}}
\def\rt'z{{\R/T'\Z}}

\def\txd{{\text{\rm d}}}

\def\Inf{{\text{\rm Inf}}}
\def\Gm{{G_\txm}}
\def\Hm{{H_\txm}}
\def\Qm{{Q_\txm}}

\def\uptwo#1{{#1^{(2)}}}

\def\up#1#2{{#1^{(#2)}}}
\def\mup#1{{^{(-1)^{#1}}}}

\def\dfn={{\overset{\text{\rm def}}\to=}}
\def\BBig(-{{\Big(\!\!-}}
\def\BBiglangle-{{\Big\langle\!\!-}}
\def\:{{\text{\rm:}}}
\def\;{{\text{\rm;}}}

\def\Linfr{{L^\infty(\R)}}

\def\bzdz\inv{{\bar z_0^{-1}}}

\nologo
\TagsOnRight
\loadeusm
\topmatter
\title{OUTER ACTIONS OF A DISCRETE AMENABLE GROUP
ON APPROXIMATELY FINITE
DIMENSIONAL FACTORS I\!I,\\
The I\!I\!I$_{\pmb\la}$-Case,\ $\pmb{\la\neq 0}$.  }
\endtitle
\rightheadtext{Outer Conjugacy I\!I}
\leftheadtext{Outer Conjugacy I\!I}
\author
{Dedicated to the Memory of Gert Kj\ae rg\aa rd  Pedersen\\
\vskip1in
Yoshikazu Katayama\\\\
Department of Mathematics,\\
Osaka Kyoiku University\\\\ and \\ \\
Masamichi Takesaki\\
\\
Department of Mathematics\\ University of California, Los Angeles}
\endauthor
\address{Department of Mathematics,
Osaka Ky$\hat\text{o}$iku University Osaka, Japan}
\endaddress
\address{Department of Mathematics,
University of California,
Los Angeles, California 90095-1555}
\endaddress

\newpage

\abstract
{To study outer actions $\a$ of a group $G$ on a factor  $\sM$ of 
type {\threel}, $0<\la<1$, we study first
the cohomology group of a group with the unitary group of an abelian 
{\vna} as a coefficient group and
establish a technique to reduce the coefficient group to the torus 
$\T$ by the Shapiro mechanism based on
the groupoid approach. We then show a functorial construction of 
outer actions of a {\cdag} on an AFD
factor of type {\threel}, sharpening the result in \cite{KtT2: \S4}. 
The periodicity of the flow of weights on
a factor $\sM$ of type {\threel} allows us to introduce an equivariant 
commutative square directly related to
the discrete core. But this makes it necessary to introduce an 
enlarged group $\Aut(\sM)_\txm$
relative to the modulus homomorphism $\txm=\mod\!\!: \Aut(\sM)\mapsto 
\rt'z$. We then discuss the
reduced {\mhjr}, which allows us to describe the invariant of outer 
action $\a$ in a simpler form than
the one for a general AFD factor: for example, the cohomology group 
$\tH_{\txm, \fs}^\out(G, N,
\T)$ of modular obstructions is a compact abelian group. Making use 
of these reductions, we prove the
classification result of outer actions of $G$ on an {\AFD} factor 
$\sM$ of type \threel.}
\endabstract
\toc\specialhead{\S0. Introduction}
\endspecialhead\specialhead{\S1. Groupoid Cohomology}
\endspecialhead
\specialhead{\S2. Model Construction I\!I}
\endspecialhead
\specialhead{\S3. Reduction of Invariants for the Case of Type
I\!I\!I$_{\la}$,
${0<\la<1.}$}
\endspecialhead
\specialhead{\S4. Outer actions of a Countable Discrete Amenable  Group on an
AFD Factor of Type I\!I\!I${_\la}$, ${0<\la<1}$.}
\endspecialhead
\specialhead{\S5. Outer actions of a Countable Discrete Amenable  Group on an
AFD Factor of Type I\!I\!I${_1}$.}
\endspecialhead
\endtoc
\endtopmatter

\document
\head{\bf \S0. Introduction}
\endhead
In this article, we continue our investigation of outer conjugacy 
classification of outer actions, say $\a$, of a {\cdg}
$G$ on a factor $\sM$ of type {\threel}, $0<\la\leq 1$. Since the 
characteristic square of the factor $\sM$ takes a
simpler form, the outer conjugacy invariants for an outer action $\a$ 
of $G$ takes a simpler form than the general case
which was completed in the last work, \cite{KtT2}. But this does not 
mean that our task was completed in the last
work. We have to reduce the general theory to the seemingly easy case 
of {\threel}, which requires more work. Once
the work is completed, we see that the final form in this particular 
case is simpler. A major hurdle for this is the fact
that the association of a discrete core $\tM_d$, a factor of type 
{\twoinf}, to a factor $\sM$ of type {\threel} is not a
functor. Accordingly, the group $\Aut(\sM)$ does not act canonically 
on $\tM_d$. The obstruction to this is the
presence of the modulus $\mod(\a)$ of $\a\in\Aut(\sM)$. Instead, an 
enlarged group $\Aut(\sM)_\txm$, which is a
central extension of $\Aut(\sM)$ by the integer group $\Z$, acts on 
$\tM_d$. The Shapiro machinery helps to
relate the characteristic square to the reduced characteristic square 
consisting of all Borel groups with compact abelian
groups, in fact, the circle group, on the crucial corner. To capture 
the Shapiro machine, we need to work with groupoid
cohomology to get a clear and natural picture, which is done in the 
first section. An interesting feature of the case of
type {\threel}, $0<\la<1$, is that the canonical two cocycle 
associated with the exact sequence:
$$\CD
0@>>>\Z@>\times T'>>\R@>>>\R/T'\Z@>>>0, \quad T'=-\log \la,
\endCD
$$
which comes from the Gauss symbol $[x], x\in \R$, the
integer $n$ such that $n\leq x< n+1$, enters naturally to the theory. 
We will use the notation
$\bracett'{\ds}$ for the cross-section:
$$
\ds=s+T'\Z\in \rtz'\mapsto \bracett'{\ds}=s-T'\left[\frac  s 
{T'}\right]\in [0, T'),\quad s\in \R.
$$

The case of type {\threeone} is even easier as its general
theory is already reduced.

Toward the completion of this article, the authors have received 
support from the Erwin Schr\"odinger
Institute and the Department of Mathematics, University of Rome, La 
Sapienza, where they visited to
work together. The authors would like to express here their gratitude 
to these institutions and
Professors K. Schmidt and S. Doplicher for their invitation and the 
hospitality extended to them, which
made this collaboration possible and enjoyable.

To keep the size of this article down, we postpone the discussion of 
examples to the subsequent paper, \cite{KtT3},
in which the third cohomology groups of easy cases are computed and 
the invariants of outer actions of such groups
are identified from their raw data.

\head{\bf \S1. Groupoid Cohomology}
\endhead
\subhead\nofrills{Shapiro's Mechanism and Dimension Shifting: }
\endsubhead
Let $\cG$ be a groupoid with $X=\up\cG0$ the space of units.
We note that whenever we consider a Borel groupoid, locally compact groupoid or
measured  groupoid, we mean by a map always a Borel map or a
measurable map. For a measured groupoid, we ignore the difference on a
null set.
By a
$\cG$-module
$\cA$, we mean a field of groups $x\in X\mapsto A(x)\in \Grp$ such that
\roster
\item"i)"  Each $g\in \cG$ gives rise to
an isomorphism $\a_g\in \Iso(A(\txs(g)), A(\txr(g))$;
\item"ii)" The family of isomorphisms $\{\a_g: g\in \cG\}$ satisfies
the chain rule:
$$
\a_{gh}=\a_g\scirc \a_h, \quad (g, h)\in \up\cG2;
$$
\item"iii)" If $x=g\in  \up\cG0$, then $\a_g=\id\in \Aut(A(x))$.
\endroster
When each $A(x), x \in X,$ is commutative, then $\cA$ is called {\it
commutative}
or {\it abelian}. We assume that $\cA$ is commutative.
An $n$-cochain, $n=0, 1, 2, \cdots,$ means a function
$$
\xi: (g_1, \cdots, g_n)\in\up\cG n\mapsto \xi(g_1, \cdots, g_n)\in
A(\txr(g_1)).
$$
The set $\tC_\a^n(\cG, \cA)$ of $n$-cochains forms a group relative to the
pointwise product. The coboundary map $\part_{n}: \tC_\a^{n}(\cG, \cA)\mapsto
\tC_\a^{n+1}(\cG,
\cA)$ is defined by
$$\aligned
(\part_{n}\xi)(g_0&, g_1, \cdots, g_n)=\a_{g_0}(\xi(g_1, \cdots, g_n))\\
&\hskip.2in\times
\prod_{k=0}^{n-1}\xi(g_0, g_1, \cdots, g_{k-2}, g_{k-1}g_k, g_{k+1}, g_{k+2},
\cdots, g_n)\mup{k-1}\\
&\hskip.5in\times
\xi(g_0, g_1, \dots, g_{n-1})\mup {n-1}\in A(\txr(g_0)).
\endaligned
$$
As usual, we have
$$
\part_n\scirc \part_{n-1}: \tZ_\a^{n-1}(\cG, \cA)\mapsto \{1\}\subset
\tZ_\a^{n+1}(\cG,
\cA).
$$
We often suppress the suffix $n$ of $\part_n$. Each element of the kernel
$\Ker(\part_n)$, denoted by $\tZ_\a^n(\cG, \cA)$, is called an
$n$-cocycle and each
element of the image $\text{\rm Im}(\part_{n-1})$, denoted by $\tB_\a^n(\cG,
\cA)$, is called an $n$-coboundary. The quotient group $\tZ_\a^n(\cG,
\cA)/\tB_\a^n(\cG, \cA)$ is called the $n$-th cohomology group of
$\cA$ and written
$\tH_\a^n(\cG, \cA)$ (See  [BN],[W1] and [W2] for the related topics).

For $n=0$, we set
$$\aligned
\tH_\a^0(\cG&, \cA)=\tZ_\a^0(\cG, \cA)\\
&=\{\xi: x\in X\mapsto \xi(x)\in A(x)\
\text{such that } \xi(\txr(g))=\a_g(\xi(\txs(g)), g \in \cG\}.
\endaligned
$$
For $n=1$,
$\tZ_\a^1(\cG, \cA)$ consists of all maps
$\xi: g\in\cG\mapsto \xi(g)\in A(\txr(g))$  such that
$$
\xi(gh)=\xi(g)\a_g(\xi(h)),\ (g, h)\in \up\cG2,
$$
and $\tB_\a^1(\cG, \cA)$ consists of all those $\xi\in\tZ_\a^1(\cG,
\cA)$ such that
$$
\xi(g)=\a_g(\eta(\txs(g))\eta(\txr(g))\inv, \quad g\in \cG,
$$
for some $\eta: x\in X\mapsto \eta(x)\in A(x)$. Each $\xi\in
\tZ_\a^1(\cG, \cA)$
gives rise to the perturbation ${}_\xi\a$ of the action $\a$ on $\cA$
given in the
following fashion:
$$
{}_\xi\a_g(u)=\xi(g)\a_g(u)\xi(g)\inv \in A(\txr(g)), \quad u\in
A(\txs(g)), g\in \cG.
$$
If $\xi\in \tB_\a^1(\cG, \cA)$, then the perturbed action ${}_\xi\a$
is conjugate to
the original action $\a$ under the group $\Int(\cA)$ of ``inner"
automorphisms of
$\cA$.

For $n=2$, each element $\xi\in \tZ_\a^2(\cG, \cA)$ is an
$\cA$-valued function on
$\up\cG2$ such that
$$
\xi(g, h)\xi(gh, k)=\a_g(\xi(h, k))\xi(g, hk)\in A(\txr(g)), \quad
(g, h, k)\in \up\cG3.
$$
The cocycle $\xi$ is a coboundary if and only if there exists a map $\eta:
g\in\cG\mapsto \eta(g)\in A(\txr(g))$ such that
$$
\xi(g, h)=\a_g(\eta(h))\eta(gh)\inv\eta(g),\quad (g, h)\in \up\cG2.
$$
Each cocycle $\xi\in \tZ_\a^2(\cG, \cA)$ gives rise to the twisted semi-direct
product groupoid:
$$
\cH=\cA\rtimes_{\a,\xi} \cG=\{(u, g)\in \cA\times \cG: u\in
A(\txr(g)), \quad g\in
\cG\}
$$
such that
$$\aligned
\up\cH2&=\{(u, g; v, h)\in \cH\times \cH: (g, h)\in \up\cG2\};\\
&\txr(u, g)=\txr(g), \quad \txs(u, g)=\txs(g);\\
(u, g)(v, h)&=(u\a_g(v)\xi(g, h), gh), \quad (u, g; v, h)\in \up\cH2.
\endaligned
$$
The original $\cG$-module $\cA$ is then viewed as a normal subgroupoid of $\cH$
and the original groupoid $\cG$ is then identified with the quotient groupoid:
$\cG=\cH/\cA$. The action $\a$ of $\cG$ on $\cA$ is then nothing but the
conjugation:
$$
\a_g(u)=(1_y, g)(u, x)(1_y, g)\inv \in A(y),\quad g=(y, g, x)\in \cG, 
$$
where $ \txr(g)=y,\ \txs(g)=x$
and $1_y$ is the identity of $A(y)$. If $\xi=\part(\eta)$, then the map:
$$
\fs_\eta: g\in \cG\mapsto (\eta(g)\inv, g)\in \cH
$$
is an injective homomorphism of the groupoid $\cG$ into $\cH$ which decomposes
$\cH$ into a semi-direct product:
$$
\cH\cong\cA\rtimes_\a \cG.
$$
For $n=0, 1, 2,$ the $\cG$-module $\cA$ does not have to be
commutative to define
$\tH_\a^n(\cG, \cA)$ as long as one is ready to give up the group
structure on the
cohomology space $\tH_\a^n(\cG, \cA), n=0, 1, 2$. For $n=2$, the
cocycle identity,
however, should be replaced by:
$$
\a_g(\xi(h, k))\xi(g, hk)\{\xi(g, h)\xi(gh, k)\}\inv=1, \quad (g, h,
k)\in \up\cG3,
$$
and the equivalence $\xi\equiv \xi'$ of two cocycles $\xi$ and $\xi'$
is defined by
the existence of $\eta\!: g\in\cG\mapsto \eta(g)\in A(\txr(g))$ such that
$$
\xi'(g, h)=\eta(g)\a_g(\eta(h))\xi(g, h)\eta(gh)\inv, \quad (g, h)\in \up\cG2.
$$

If the groupoid $\cG$ is a topological groupoid and $\cA$ admits a topological
structure such that all the operations are continuous, then we
request that cocycles
are all Borel. To demand the continuity on cocycles is too
restrictive as seen in the
group case. If $\cG$ is a measured groupoid, then all the identities
mean to hold
almost everywhere relative to the relevant measure class.

\proclaim{Proposition 1.1} Let $\cA$ be a $\cG$-module. For each
$y\in X=\up\cG0$,
let $B(x)$ be the set of all $A(x)$-valued functions on
$\cG_x=\txs\inv(x)$ and set
$$
\cB=\bigcup_{x\in X}^\bullet B(x)=\{b: \cG\mapsto \cA\ \text{such
that } b(g)\in
A(x), g=(y, g, x) \in \cG\}.
$$
For each $g=(y, g, x)\in\cG$, define the map $\b_g\!\!: B(x)\mapsto
B(y)$ in the
following fashion{\rm:}
$$\aligned
(\b_gb)(h)=\a_g(b(hg))\in A(y), \quad b\in B(x), (h, g)\in \up\cG2.
\endaligned
$$
Then for each $n=2, 3, \cdots,$
$$
\tH_\b^n(\cG, \cB)=\{1\}.
$$
More explicitly, if $\xi\in \tZ_\b^n(\cG, \cB)$, then $\eta\in
\tC_\b^{n-1}(\cG, \cB)$
defined by
$$\aligned
\eta(g&; g_1, g_2, \cdots, g_{n-1})\\
&=\ainv{g}(\xi(\txr(g); g, g_1, \cdots, g_{n-1}))\in
A(\txs (g)), \quad (g, g_1,
\cdots, g_{n-1})\in \up\cG n.
\endaligned
$$
gives
$$
\xi=\part \eta.
$$
\endproclaim
\demo{Proof} The cocycle identity:
$$\aligned
1&=\b_{g_0}(\xi(g_1, \cdots, g_n))\prod_{i=1}^n\xi(g_0, g_1, g_2,
\cdots, g_{i-2},
g_{i-1}g_i, g_{i+1}, \cdots, g_n)\mup i\\
&\hskip 1in\times
\xi(g_0, \cdots, g_{n-1})\mup{n}, \quad (g_0, g_1, \cdots, g_n)\in \up\cG{n+1}
\endaligned
$$
gives
$$\aligned
\xi&(g_1, \cdots, g_n)=\b_{g_0}^{-1}\Big(\prod_{i=1}^n\xi(g_0, g_1,
g_2, \cdots,
g_{i-2}, g_{i-1}g_i, g_{i+1}, \cdots g_n)\mup {i-1}\\
&\hskip 1.5in\times
\xi(g_0, \cdots, g_{n-1})\mup{n-1}\Big)
\endaligned
$$
which means that for each $(g, g_1, \cdots, g_n)\in\up\cG {n+1}$ with
  $g=g_0$,
$$\aligned
\xi(g&; g_1, \cdots, g_n)\\
&=\ainv{g_0}\Big(\prod_{i=1}^n\xi(gg_0^{-1}; g_0, g_1,
g_2,
\cdots, g_{i-2}, g_{i-1}g_i, g_{i+1}, \cdots g_n)\mup {i-1}\\
&\hskip 1.5in\times
\xi(gg_0^{-1}; g_0, \cdots, g_{n-1})\mup{n}\Big)\\
&=\ainv{g}\Big(\prod_{i=1}^n\xi(\txr(g); g, g_1, g_2,
\cdots, g_{i-2}, g_{i-1}g_i, g_{i+1}, \cdots g_n)\mup {i-1}\\
&\hskip 1.5in\times
\xi(\txr(g); g, g_1, \cdots, g_{n-1})\mup{n}\Big).\\
\endaligned
$$
Setting
$$\aligned
\eta&(g; g_1, \cdots, g_{n-1})=\ainv{g}\Big(\xi(\txr(g); g, g_1,
\cdots, g_{n-1})\Big)
\in A(\txs(g));\\
&\hskip 1in\eta\in\tC_\b^{n-1}(\cG, \cB),
\endaligned
$$
we compute for $(g, g_1, \cdots, g_n)\in \up\cG{n+1}$ with $y=\txr(g)$
$$\allowdisplaybreaks\aligned
(\part \eta)&(g; g_1, \cdots, g_n)=\a_{g_1}\Big(\eta(gg_1; g_2,
\cdots, g_n)\Big)\\
&\hskip .2in\times
\prod_{i=1}^{n-1}\eta(g; g_1, \cdots, g_{i-1}, g_ig_{i+1}, g_{i+2},
\cdots, g_n)\mup i\\
&\hskip .5in\times
\eta(g; g_1, g_2, \cdots, g_{n-1})\mup{n}\\
&=\a_{g_1}\Big(\ainv{gg_1}\big(\xi(y; gg_1, g_2, \cdots, g_n)\big)\Big)\\
&\hskip .2in\times
\prod_{i=1}^{n-1}\ainv g\big(\xi(y; g, g_1, \cdots, g_{i-1},
g_ig_{i+1}, g_{i+2},
\cdots, g_n)\big)\mup i\\
&\hskip .5in\times
\ainv g\big(\xi(y; g, g_1, g_2, \cdots, g_{n-1})\big)\mup{n}\\
&=\xi(g; g_1, \cdots, g_n).
\endaligned
$$
This completes the proof.
\qed
\enddemo

Each $A(x)$ is a submodule of $B(x)$ for each $x\in X$, hence we get an exact
sequence:
$$\CD
\{1\}_x@>>>A(x)@>i_x>>B(x)@>j_x>\underset{{\fs\!_j}_x}\to\longleftarrow>C(x)
@>>>\{1\}_x\qquad (x\in X),
\endCD
$$
and another $\cG$-module $\cC$:
$$
\cC=\bigcup_{x\in X}^\bullet C(x), \quad \g_g(bA(x))\overset{\text{def}}\to{=}\b_g(b)A(y),
b\in B(x), g=(y, g,
x)\in \cG.
$$
Symbolically we can write the exact sequence of $\cG$-modules:
$$\CD
X@>>>\cA@>i>>\cB@>j>>\cC@>>>X.
\endCD
$$
Take $\xi\in \tZ_\a^n(\cG, \cA)$. With
$$
(i_*\xi)(g; g_1, \cdots, g_n)=\xi(g_1, \cdots, g_n)\in A(\txr(g_1)),
\quad (g, g_1, \cdots,
g_n)\in \up\cG {n+1},
$$
we obtain a cocycle $i_*\xi\in \tZ_\b^n(\cG, \cB)=\tB_\b^n(\cG, \cB)$
by Proposition
1.1. Thus there exists a cochain $\eta\in \tC_\b^{n-1}(\cG, \cB)$ such that
$$
i_*\xi=\part \eta.
$$
Then set
$$\aligned
\z(x; g_1, \cdots, g_{n-1})&=j_x(\eta(x; g_1, \cdots, g_{n-1}))\\
&=\eta(x; g_1, \cdots, g_{n-1})A(x)\in C(x).
\endaligned
$$
Since $\part j=j\part$, we get $\z\in \tZ_\g^{n-1}(\cG, \cC)$ and naturally
$$
\tH_\a^n(\cG, \cA)\cong \tH_\g^{n-1}(\cG, \cC)
$$
under the map: $[\xi]\in\tH_\a^n(\cG, \cA)\mapsto
[\z]\in\tH_\g^{n-1}(\cG, \cC)$.
Summarizing the above discussion, we obtain:

\proclaim{Proposition 1.2 (Dimension Shifting)} If $\{\cA,\a\}$ is a
$\cG$-module, then
there exists a natural $\cG$-module $\{\cC, \g\}$ and a natural
isomorphism{\rm:}
$$
\tH_\a^n(\cG, \cA)\cong \tH_\g^{n-1}(\cG, \cC).
$$
\endproclaim

\subhead{Pullback, Reduction and Induction}
\endsubhead
Let $\cH$ be a groupoid with $Y=\up\cH0$. Suppose that $f$ is a
surjective map from a space
$X$ onto $Y$ (if applicable, we assume that the map $f$ is Borel).
Then we have a fibration of $X$:

\hsize= 2.75in
$$
X=\bigcup_{y\in Y}^\bullet X(y), \ X(y)=f\inv (y), y \in Y.
$$
Then we set
$$\allowdisplaybreaks\aligned
\cG=\bigcup_{(y, z)\in Y\times Y}\{X(z)\times\cH_z^y\times X(y)\},
\endaligned
$$
where $\cH_y^z=\txr\inv(z)\cap \txs\inv(y), (z, y)\in Y^2$. We then define the
range and the source maps and the product in $\cG$ as follows:

\hsize=5.5 true in

\vskip.1in
\hskip2.75in\special{picture Fig1 scaled .5}

\noindent
$$\allowdisplaybreaks\aligned
&\hskip.5in \txr(z, h, x)=z,\quad \txs(z, h, x)=x;\\
&(z, g, y)(y, h, x)=(z, gh, x),\quad (g, h)\in\up\cH2,\\
& z\in X(\txr(g)), y\in X(\txs(g))=X(\txr(h)),  x\in X(\txs(h)).
\endaligned
$$

{\smc Definition 1.3.} The groupoid $\cG$ is called the {\it pullback} of $\cH$ by the map
$f$ and denoted by $f^*(\cH)$.

The map:
$$
f_*: (z, h, x) \in \cG\mapsto f_*(z, h, x)=(f(z), h, f(x))\in \cH
$$
is a groupoid homomorphism of $\cG$ onto $\cH$.

If $f_i, i=1, 2,$ are maps from $X_i$ onto $Y$, then we have the
fiber product $X=X_1*X_2$
relative to $f_1$ and $f_2$:
$$
X=\{(x_1, x_2)\in X_1\times X_2: f_1(x_1)=f_2(x_2)\},
$$
and the map $f: x=(x_1, x_2)\in X\mapsto f(x)=f_1(x_1)=f_2(x_2)\in Y$
which makes the
following diagram commutative:
$$\CD
X@>\pr_1>>X_1\\
@V\pr_2VV@V f_1VV\\
X_2@>f_2>>Y
\endCD\qquad f=f_1\scirc \pr_1=f_2\scirc \pr_2.
$$
The pullbacks $\cG=f^*(\cH)$, $\cG_1=f_1^*(\cH)$, $\cG_2=f_2^*(\cH)$ and
$\cH$ form
the commutative diagram:
$$\CD
\cG@>\pr_{1*}>>\cG_1\\
@V\pr_{2*}VV@V f_{1*}VV\\
\cG_2@>f_{2*}>>\cH
\endCD\qquad f_*=f_{1*}\scirc \pr_{1*}=f_{2*}\scirc \pr_{2*}.
$$

Let $Y\subset X$ be a subset  such that the saturation
$[Y]=X$, i.e.,
$$
[Y]=\txs\scirc \txr\inv(Y)=X,
$$
equivalently for every $x\in X$ there exists $g=(y, g, x)\in \cG$ such that
$y\in Y$. If applicable, we assume that
the set $Y$ is a Borel
subset of $X$.

{\smc Definition 1.4.} In the above setting, let $\cG_Y$ be the set
of all those $g\in \cG$
such that $\txr(g)\in Y$ and $\txs(g)\in Y$, i.e.,
$$
\cG_Y=\txr\inv(Y)\cap \txs\inv(Y).
$$

\proclaim{Proposition 1.5} If $Y$ is a subset of $X$ such that
$X=[Y]$, then there exists a
surjective map $f: X\mapsto Y$ such that $\cG$ is naturally
isomorphic to the pullback
groupoid $f^*(\cG_Y)$.
\endproclaim
\demo{Proof}
Consider the map
$$
\txs\!: g\in \txr\inv(Y)\mapsto x=\txs(g)\in X
$$
and its cross-section (if applicable, we take a Borel
cross-section $\g$)
$$
\g: x\in X\mapsto \g(x)\in \txr\inv(Y),
$$
such that $\g(y)=y$ if $y\in Y$.
Set
$$\aligned
f(x)&=\txr(\g(x))\in Y, \quad x\in X.
\endaligned
$$
We claim $\cG\cong f^*(\cG_Y)$. For each $(z, g, x)\in f^*(\cG_Y)$, set
$$
\pi(z, g, x)=\g(z)\inv g\g(x)\in \cG
$$
which makes sense because
$$\allowdisplaybreaks\aligned
\txs(\g(z)\inv)&=\txr(\g(z))=f(z)=\txr(g)\in Y;\\
\txr(\g(x))&=f(x)=\txs(g).
\endaligned
$$
For each $((z, g, y), (y, h, x))\in \up{f^*(\cG_Y)}2$, we get
$$\allowdisplaybreaks\aligned
\pi((z, g, y)(y, h, x))&=\pi(z, gh, x)=\g(z)\inv gh\g(x)\\
&=\g(z)\inv g\g(y)\g(y)\inv h\g(x)\\
&=\pi(z, g, y)\pi(y, h, x).
\endaligned
$$
Hence the map $\pi$ is multiplicative. The inverse $\pi\inv$ is given by:
$$\allowdisplaybreaks\aligned
\pi\inv(z, g, x)&=(z, \g(z)g \g(x)\inv, x), \quad (z, g, x)\in \cG.
\endaligned
$$
This proves the assertion.
\qed
\enddemo

Similarly, if $\cB$ is an $\cH$-module, (not necessarily
commutative), then the
surjective map $f \ \coron X\mapsto Y$ also gives rise to the pullback
$f^*(\cH)$-module
$\cA=f^*(\cB)$ in the following way:
$$\allowdisplaybreaks\aligned
\CD A(x)@>\a_{(z, g, x)}>>A(z)\\
@|@|\\
B(f(x))@>\b_g>> B(f(z))
\endCD \qquad (z, g, x)\in
\cG=f^*(\cH),
\endaligned
$$
where $\b$ is the action of $\cH$ on $\cB$.

{\smc Definition 1.6.} i) If $\cH=\cG_Y, Y\subset X,$ and the map $f$ is
given by Proposition  1.5, then the  above $\cG$-module $f^*(\cB)$ is
called the {\it
induced} $\cG$-module and written $\cA=\Ind_Y^X\cB$ or $\cA=\Ind_{Y\uparrow
X}\cB$.

ii) If $\cA$ is a $\cG$-module, then
$$
\cA_Y=\bigcup_{y\in Y}^\bullet A(y)
$$
is naturally $\cG_Y$-module, which will be called the {\it reduced} module over
$\cG_Y$ or the {\it reduced} $\cG_Y$-module.

\proclaim{Proposition 1.7} If $X=[Y]$, then every $\cG$-module $\cA$,
not necessarily
commutative,  is obtained from the reduced
$\cG_Y$-module $\cA_Y$ as the induced  module.
\endproclaim
The proof is exactly the same as  Proposition 1.5, and we
leave details to the
interested reader.

\proclaim{Proposition 1.8} If $Y\subset X$ and $X=[Y]$, then the embedding map
$$
i_Y\coron \cG_Y\hookrightarrow \cG
$$
gives  rise to the pullback
map, i.e.,  the restriction map, with the following properties{\rm:}
$$\aligned
&i_Y^*: \xi\in\tZ_\a^n(\cG, \cA)\mapsto \xi_Y=\xi|_Y\in
\tZ_\a^n(\cG_Y, \cA_Y);\\
f_*^*&\scirc i_Y^*(\xi)\equiv \xi\quad \mod\ \tB_\a^n(\cG, \cA), \quad \xi\in
\tZ_\a^n(\cG, \cA);\\
i_Y^*\scirc f_*^*&(\xi)\equiv \xi\quad \mod\ \tB_\a^n(\cG_Y, \cA_Y),
\quad \xi\in
\tZ_\a^n(\cG_Y, \cA_Y);\\
&\hskip .5in i_Y^*\coron \tH_\a^n(\cG, \cA)\cong \tH_\a^n(\cG_Y, \cA_Y),
\endaligned\nopagebreak
$$
where  $f$ is the map of Proposition 1.5 and  the map $f_*: \cG\longmapsto \cG_Y$ is 
given by
$f_*(g)=\gamma (z) g \gamma (x)^{-1},\ g=(z,g,x)\in\cG$.
\endproclaim
\demo{Proof} In view of the last two propositions, we may and do
assume that the
groupoid $\cG$ and the
$\cG$-module $\cA$ are both obtained as $\cG=f^*(\cG_Y)$  and
$\cA=f^*(\cA_Y)$. In this setting, we get $A(x)=A(f(x))$ and
$\a_{\g(x)}=\id_{A(x)}$,
the  identity map on $A(x)$, for every $x\in X$.

For $n=0$, each $\xi\in \tZ_\a^0(\cG, \cA)$ is a map $\xi: x\in X\mapsto
\xi(x)\in A(x)$ such that
$$
\xi(\txr(g))=\a_g(\xi(\txs(g))), \quad g\in \cG.
$$
The restriction $\xi_Y$ satisfies the same identities for $\cG_Y$ and
$\xi(x)=\xi_Y(f(x)),
\allowmathbreak x \in X$. Hence $\xi_Y\in \tZ_\a^0(\cG_Y, \cA_Y)$.

If $\eta\in  \tZ_\a^0(\cG_Y, \cA_Y)$,  then $\xi= f^*_*\eta$ satisfies
for each $g=(z, f_*(g),
x)\in \cG$, we have
$$
g=\g(z)\inv f_*(g)\g(x),
$$
and
$$\allowdisplaybreaks\aligned
\xi(z)&=\eta(\txr(f_*(g)))
=\a_{f_*(g)}\Big(\eta(f(x))\Big)\\
&=\a_{f_*(g)}(\xi(x))=\a_g(\xi(x)).\\
\endaligned
$$
Hence we get $\xi= f^*_*\eta\in \tZ_\a^0(\cG, \cA)$.

Since $\xi(x)=\xi(f(x)), x \in X,$ for
every $\xi\in \tZ_\a^0(\cG, \cA)$, we conclude that
   $\tH_\a^0(\cG, \cA)\cong \tH_\a^0(\cG_Y, \cA_Y)$ under the
isomorphism $i_Y^*$.

The case $n=1$: Each $\xi\in \tZ_\a^1(\cG, \cA)$ satisfies
$$
\xi(gh)=\xi(g)\a_g(\xi(h)), \quad (g, h)\in \up\cG2.
$$
The restriction $\xi_Y$ satisfies the same identity, so that it is a cocycle in
$\tZ_\a^1(\cG_Y, \cA_Y)$. Now choose $\xi_Y\in \tZ_\a^1(\cG_Y, \cA_Y)$ and set
$$
\xi(z, g, x)=(f_*^*\xi_Y)(z, g, x)=\xi_Y(\g(z) g\g(x)\inv).
$$
For each pair $g=(z, g, y), h=(y, h, x)\in \cG$, we have
$$\allowdisplaybreaks\aligned
\xi(gh)&=\xi_Y(f_*(gh))=\xi_Y(f_*(g)f_*(h))\\
&=\xi_Y(f_*(g))\b_{f_*(g)}(\xi_Y(f_*(h))=\xi(g)\a_g(\xi(h)).
\endaligned
$$
Hence $f_*^*(\xi_Y)\in \tZ_\a^1(\cG, \cA)$.

Suppose $\xi_Y=i_Y^*(\xi)$. We then compare $\xi$ and $f_*^*\xi_Y$.
For $g=(z, g,
x)$, we write $f_*(g)=(f(z), f_*(g), f(x))=\g(z)g\g(x)\inv\in \cG_Y$
and compute:
$$\allowdisplaybreaks\aligned
(f_*^*\xi)(g)&=\xi(\g(z)g\g(x)\inv)=\xi(\g(z))\xi(g\g(x)\inv)\quad(\text{as }
\a_{\g(z)}=\id)\\
&=\xi(\g(z))\xi(g)\a_g\big(\xi(\g(x)\inv)\big)\\
&=\xi(\g(z))\xi(g)\a_{g}\big(\ainv{\g(x)}(\xi(\g(x))\inv)\big)\\
&=\xi(\g(z))\xi(g)\a_{g}\Big(\xi(\g(x))\inv\Big).
\endaligned
$$
Hence the above calculation becomes:
$$\allowdisplaybreaks\aligned
(f_*^*\xi)(g)&=\xi(\g(z))\xi(g)\a_g(\xi(\g(x))\inv).
\endaligned
$$
Therefore we get $\xi\equiv f_*^*\xi\quad \mod\ \tB_\a^1(\cG, \cA)$.

For $n=2, 3, \cdots,$ we use the dimension shifting, Proposition 1.2. From the
construction of $\{\cC, \g\}$ from $\{\cA, \a\}$, it follows that the reduced
$\cG_Y$-modules $\cC_Y$ and $\cA_Y$ are  related in the exactly same way as the
original modules $\cA$ and $\cC$ are. Hence we have by mathematical induction:
$$
\tH_\a^n(\cG, \cA)\cong \tH_\g^{n-1}(\cG, \cC)\cong\tH_\g^{n-1}(\cG_Y, \cC_Y)
\cong \tH_\a^n(\cG_Y, \cA_Y).
$$
Tracing the isomorphims, we conclude that the isomorphism is indeed
given by $f_*^*$
and $i_Y^*$.
\qed
\enddemo
{\smc Remark 1.9.} In the case that $Y$ is a singleton set $\{y_0\}$,
the reduced
groupoid
$\cG_Y$ is a group,  say $H$.  The associated principal groupoid
$\widetilde\cG$, the
equivalence relation groupoid given by the orbit structure of $\cG$,
is transitive. This is
precisely the conventional induction procedure and also the Shapiro
mechanism. This case
is also relevant to us in the case of a system based on a factor of
type {\threel} as will be
seen in the later sections.

{\smc Definition 1.10.} A {\it normal subgroupoid} $\cN$ of $\cG$ is a field
$$
\{N(x)\i \txr\inv(x)\cap\txs\inv(x): x \in X\}
$$
of groups such that
$$
N(y)=gN(x)g\inv, \quad (y, g, x)\in \cG.
$$

For a commutative $\cG$-module $\cA$ with trivial action of $\cN$, we
can define the
group $\tZ_\a(\cG, \cN, \cA)$ of characteristic cocycles $(\la, \mu)$
as  in \cite{ST1} and the
group $\La_\a(\cG, \cN,
\cA)$ of characteristic invariants. Each $(\la, \mu)\in \tZ_\a(\cG,
\cN, \cA)$ gives rise to an
exact sequence of $\cG$-modules equipped with cross-section $\fs$:
$$\CD
\cE\coron X@>>>\cA@>i>>\cE=\cA\rtimes_{\a, \mu}
\cN@>j>\underset{\fs}\to\longleftarrow >\cN@>>>X
\endCD
$$
where
$$\allowdisplaybreaks\aligned
&\fs: (m, x) \in N(x)\mapsto \fs(m, x)=(1_x,
m)\in\cE(x)=A(x)\rtimes_{\a, \mu_x} N(x);\\
&\qquad\fs(m, x)\fs(n, x)=\mu(m, n; x)\fs(mn,  x)=\mu_x(m, n)\fs(mn,  x);\\
& \a_{z, g, x}(\fs(g\inv mg,  x))=\la(m; z, g, x)\fs(m ,z), \quad (z, g, x)\in
\cG, (m, x)\in N(x).
\endaligned
$$
Conversely, if we have an exact sequence of $\cG$-modules:
$$\CD
\cE\coron\ X@>>>\cA@>i>>\cE@>j>>\cN@>>>X,
\endCD
$$
then a cross-section $\fs$ (if applicable,  we take a Borel
cross-section) of the map $j$  gives rise to $\cA$-valued
functions $\la$ on
$\cN*_\txr\cG=\{(m,  g)\in
\cN\times  \cG: \txs(m)=\txr(g)\}$ and $\mu$ on $\up\cN2$ such that
$$\allowdisplaybreaks\aligned
&\fs(m, x)\fs(n, x)=\mu(m, n; x)\fs(mn, x), \quad (m, x), (n, x)\in N(x);\\
&\a_{z, g, x}(\fs(g\inv m g, x))=\la(m; z, g, x)\fs(m, z),  ((m, x),
(z, g, x))\in \cN_x\times
\cG_x^z,
\endaligned
$$
and the pair $(\la,  \mu)$ falls in the group $\tZ_\a(\cG, \cN,
\cA)$. We denote the
collection of the conjugate classes of exact sequences by $\X(\cG,
\cN, \cA)$ and each
exact sequence a {\it crossed extension of $\cA$ by $\cN$}.
   The group
multiplication in $\tZ_\a(\cG, \cN, \cA)$ reflects to the following
operations in $\X(\cG,
\cN, \cA)$:
\roster
\item"i)" For any two crossed extensions $\cE_1, \cE_2\in \X(\cG, \cN, \cA)$:
$$\CD
\cE_1\coron X@>>>\cA@>i_1>>\cE_1@>j_1>>\cN@>>>X;\\
\cE_2\coron X@>>>\cA@>i_2>>\cE_2@>j_2>>\cN@>>>X
\endCD
$$
the product crossed extension $\cE$ is defined to be the quotient
module  of the fiber
product:
$$\aligned
\cE&=\cE_1*\cE_2\\
&=\{(e_1, e_2)\in \cE_1\times \cE_2: j_1(e_1)=j_2(e_2)\}/\{(i_1(a),
i_2(a\inv)): a\in \cA\}.
\endaligned
$$
\item"ii)" The inverse $\cE\inv$ is then given by:
$$\CD
\cE'\coron\ X@>>>\cA@>\{a\in \cA\mapsto i(a\inv)\}>>\cE@>>>\cN@>>>X
\endCD
$$

\endroster

\proclaim{Proposition 1.11} If $Y\i X$ is a subset of $X$ such that
$X=[Y]=\txs(\txr\inv(Y))$, then with the map  $f$ of {\rm Proposition
1.5} the maps
$f_*^*$ and
$i_Y^*$ give isomorphisms between $\La_\a(\cG, \cN, \cA)$ and
$\La_\a(\cG_Y,\allowmathbreak\cN_Y, \cA_Y)$ which are the inverse of
one another.
\endproclaim
\demo{Proof} First, we may and do assume that $\cA=\Ind_Y^X \cB$ and
$\cG=f^*(\cH)$
with $\cB=\cA_Y$ and $\cH=\cG_Y$. Also the normal subgroupoid $\cN$ is a
$\cG$-module and therefore it is conjugate to the  induced
$\cG$-module $\Int_Y^X (\Cal
M)$ with $\Cal  M=\cN_Y$ as a
$\cG$-module. However, it is not entirely trivial to relate the
structure of the inclusions:
$\Cal M\subset \cG_Y$ and $\cN\subset \cG$. We have to study the way
that the entire
groupoid $\cG$ is related to the reduced one $\cG_Y$. As
$\cG=f^*(\cG_Y)$, we have
$$\aligned
&\quad\cG_x^z=\{z\}\times \cH_{f(x)}^{f(z)}\times \{x\};\\
   \cG_x^x&=\{x\} \times \cH_{f(x)}^{f(x)}\times
\{x\}=\g(x)\inv\cH_{f(x)}^{f(x)}\g(x);\\
N(x)&=\g(x)\inv N(f(x))\g(x)=\{x\}\times M(f(x))\times\{x\}.
\endaligned
$$
Hence we have $\cN=f^*(\cM)$.

Now every $\cE\in \X_\a(\cG, \cN, \cA)$ conjugate to a crossed
extension of the form
$f^*(\cF)$ with $\cF\in \X_\a(\cH, \cM, \cB)$ and $\cF$ is given as
$\cF=\cE_Y$. It is
straightforward now to see that $f^*(\cF_1)\cong f^*(\cF_2)$ if and
only if $\cF_1\cong
\cF_2$. Hence the $f_*^*$ and $i_Y^*$ are isomorphisms between
$\La_\a(\cG, \cN,
\cA)$ and $\La_\a(\cG_Y, \cN_Y, \allowmathbreak\cA_Y)$ which are
inverse of one another.
\qed
\enddemo

\subhead\nofrills{\bf Non-Polish Groupoid:}
\endsubhead\quad
In many cases, we encounter the following situation:
\roster
\item"i)" An ergodic flow $\{\sA, \R, \th\}$ is given on an
abelian separable
{\vna } $\sA$. Let $\{X, \txm\}$ be the measure theoretic spectrum of $\sA$,
i.e., $\{X,
\txm\}$ is a standard measure space such that $\sA=L^\infty(X, \txm)$;
\item"ii)" A Polish group $H$ acts on the flow $\{\sA, \R, \th\}$ via
$\a$, i.e., $\a$ is
a homomorphism of $H$ into the group $\Aut_\th(\sA)=\{\sig \in \Aut(\sA):
\sig\scirc \th_s=\th_s\scirc \sig, s \in \R\}$, which gives rise to a
joint action of
$\wtH=H\times \R$ on $\sA$. We denote it by $\a$ also and by $\th$ if
we restrict
$\a$ to the action of $\R$. We also use the notations:
$$
(\a_{g, s}f)(x)=f(T_{g, s}^{-1}x), \quad f\in \sA=L^\infty(X, \txm),
(g, s)\in \wtH.
$$
\item"iii)" A normal subgroup $L$ of $H$ contained in the kernel
$L\subset \Ker(\a)$
is given, so that the action of $H$ factors through the quotient group $Q=H/L$;
\item"iv)" The normal subgroup $L$ is {\it not} closed, so that the
quotient group
$Q$ does not have a reasonable topological or Borel structure despite its
significance in the theory.
\endroster
The action $\a$ of $\wtH$ on $\sA$ gives rise to the action of
$\wtQ=Q\times \R$
which will be denoted by $\a$ again. Now the groupoid $\cG_\wtQ=X\rtimes \wtQ$
is the groupoid which is relevant to us despite the lack of a
reasonable Borel structure
on it. In what follows, we consider the discrete topology on $Q$ side
and the usual
topological and Borel structures on $\R$-side. Namely we are going to
consider those
functions $f$ on $\cG_\wtQ$ such that the map: $(x, s)\mapsto f(x, q,
s)$ is jointly
Borel as a function on $X\times \R$ for each fixed $q\in Q$. Namely,
we consider
the product Borel structure of those on $X$ and $\R$ and the discrete
Borel structure
on $Q$. A typical example will be the automorphism group $\Aut(\sM)$  of  a
separable factor $\sM$ as $H$ and $\cntr(\sM)$ as $L$ and the flow $\{\sC, \R,
\th\}$ of weights on $\sM$ as $\{\sA, \R, \th\}$.

\head{\bf \S2. Model Construction I\!I}
\endhead
Let $\{\sC, \R, \th\}$ be an ergodic flow to be fixed throughout this
section. Let $H$ be a {\cdg} and $\a$ an action of $H$ on the flow
$\{\sC, \R, \th\}$, i.e., a homomorphism of $H$ into the group
$\Aut_\th(\sC)$ of all automorphisms of $\sC$ commuting with the
flow $\th$. Let $\{X, \txm\}$ be the measure spectrum of $\sC$, i.e., a
standard $\sig$-finite measure space on which $\R$ acts as a
one-parameter group of non-singular transformations $\{T_s: s \in
\R\}$, so that
$$
\th_s f(x)=f(T_s^{-1}x), \quad f\in \sC=L^\infty(X, \txm), x\in X, s \in
\R.
$$
We denote the action $\a$ of $H$ on the space $X$ in the following
fashion:
$$
\a_g f(x)=f(xg), \quad g\in H, f \in \sC.
$$
When  we consider the joint action of $\wtH=H\times \R$ on $X$,
denoted simply by $\a$ also, we write
$$
(\a_{\tilde g}f)(x)=f(x\tilde g)=f(T_s^{-1}xg), \quad \tilde g=(g, s)\in
\wtH=H\times\R.
$$
Let $L$ be a normal subgroup of $H$ contained in the Kernel
$\Ker(\a)$ of $\a$, i.e., $L$ does not act on $\sC$ at all.

This section will be devoted to a construction of  an action $\a$ of $H$
on a separable strongly stable factor $\sM$ for any characteristic
invariant $\chi\in \La_\a(\wtH, L, \sU(\sC))$ such that
\roster
\item"i)" The flow of weights on $\sM$ is conjugate to the flow $\{\sC,
\R, \th\}$, which will be identified;
\item"ii)" The modulus  $\mod(\a)$ is precisely the preassigned action
$\a$ on $\sC$ of $H$;
\item"iii)"  $$L=\a\inv(\cntr(\sM));$$
\item"iv)" $$\chi=\chi(\a).$$
\endroster
Here the strong stability of $\sM$ means that $\sM\cong\sM\botimes
\sR_0$ with $\sR_0$
an approximately finite dimensional factor of type {\twoone}. This is
equivalent to the non
commutativity of the quotient group $\Intb(\sM)/\Int(\sM)$ of the
group $\Intb(\sM)$ of
approximately inner automorphisms by the group $\Int(\sM)$ of inner
automorphisms.

The  joint action $\wtH$ on $X$ gives rise to a standard measured
groupoid $\Cal H=X\times \wtH$ such that
$$\aligned
&\txr(y, \tilde g)=y\quad \txs(y, \tilde g)=y\tilde g, \quad y\in X, \tilde
g=(g, s)\in \wtH=H\times \R;\\
&\qquad (y, \tilde g)(y\tilde g, \tilde h)=(y, \tilde g\tilde h),\quad \tilde
h=(h, t)\in \wtH.
\endaligned
$$
In order to shorten notations we write
$$\tilde g\in \cG\text{ and } \tilde h\in
\cG\text{ for } (y, \tilde g)\in \cG\text { and }(x, \tilde h)\in \cG
$$
respectively
omitting  the range $y=\txr(\tilde g), x=\txr(\tilde h)\in X$. At the same
time, to specify the range and the source explicitly, we often write
$$
(y, \tilde g, x)\in \cG\text{ for } y=\txr(y, \tilde g, x) \text{ and
} x=\txs(y,
\tilde g, x),\text{ i.e., } x=y\tilde g.
$$

For each $x\in X$, let $K(x)$ be the isotropy group of $x$, i.e.,
$$
K(x)=\{\tilde g\in \wtH: x=x\tilde g\}.
$$
The map: $x\in X\mapsto K(x)\in  \Sub(\wtH)$ is a Borel map from the
standard Borel space $X$ into the standard Borel space $\Sub(\wtH)$ of all
closed subgroups of $\wtH$ such that $K(x)=K(T_sx), s\in \R, x \in X,$
since the flow $\th$ and the joint action $\a$ of $\wtH$ commute. The
ergodicity of
$\th$ then implies that $K(x)=K\in  \Sub(\wtH)$ for some fixed closed
subgroup $K$  of $\wtH$. Then of course, $K$ is the Kernel $\Ker(\a)$  of
the joint action
$\a$ of $\wtH$, hence it is a normal closed subgroup of $\wtH$ which
contains the normal subgroup $L$. Let $\cN=X\times L$ denote  the
normal  subgroupoid of $\cG$.

\proclaim{Proposition 2.1} Let $A$ be the unitary group $\sU(\sC)$ of $\sC$.
Then there is a natural isomorphism$:(\la, \mu)\in \tZ_\a(\wtH, L,
A)\mapsto (\tla,
\tmu)\in \tZ(\cH, \cN, \T)$ from the group $\tZ_\a(\wtH, L, A)$ of
characteristic
cocycles onto the group
$\tZ(\cG, \cN, \T)$ of characteristic cocycles which maps precisely
the group $\tB_\a(\wtH, L, A)$ of coboundaries onto the subgroup
$\tB(\cG, \cN, \T)$ of coboundaries, so that it induces a natural
isomorphism{\rm:}
$$
\chi\in\La_\a(\wtH, L, A)\mapsto \tchi\in\La(\cG, \cN, \T).
$$
\endproclaim
\demo{Proof} To each $(\la, \mu)\in \tZ_\a(\wtH, L, A)$, we want to
associate a cocycle $(\tilde \la, \tilde \mu)\in \tZ(\cG, \cN, \T)$. First,
we realize $(\la(m, \tilde g), \mu(m, n))$ for $m, n \in L$ and $\tilde
g=(g, s)\in \wtH$ as  $\T$-valued Borel functions over $X$ so that the
cocycle identities hold almost everywhere and write them:
$$
\la(x; m; \tilde g)\quad\text{and}\quad \mu(x; m, n),\quad x\in X,\ m, n
\in L,\ \text{and}\ \tilde g=(g, s)\in \wtH.
$$
Then set
$$\aligned
&\tilde \mu((x, m); (x, n))=\mu(x; m, n), \quad ((x, m); (x, n))\in
\cN^{(2)};\\
&\tilde \la((x, m); (x, \tilde g))=\la(x; m; \tilde g), \quad ((x, m); (x,
\tilde g))\in \cG^{(2)}.
\endaligned
$$
Let us check the cocycle identity:
\roster
\item"i)" The 2-cocycle identity for $\mu$ guarantees the 2-cocycle
identity for $\tilde \mu$ almost everywhere;
\item"ii)" If $m=(z, m, z)\in \cN$ and $(\tilde g=(z, \tilde g, y),
\tilde h=(y,
\tilde h, x))\in \uptwo{\cG}$, then $y=z\tilde g, \ x=y\tilde h$ and
$$\aligned
&\tilde\la( m; \tilde g\tilde h)=\tilde \la((z, m, z); (z, \tilde g\tilde
h, x))=\la(z; m; \tilde g\tilde h)\\
&\hskip .5in=\la(z; m; \tilde g)\la(z\tilde g; \tilde g\inv m\tilde
g; \tilde h)\\
&\hskip .5in=\tilde \la((z, m, z); (z, \tilde g, z\tilde
g))\tilde\la((z\tilde g,
\tilde g\inv m \tilde g, z\tilde g); (z\tilde g, \tilde h, x))\\
&\hskip.5in=\tilde \la((z, m, z); (z, \tilde g, y))\tilde \la((y, \tilde g\inv
m\tilde g, y); (y, \tilde h, x))\\
&\hskip .5in=\tilde \la(m; \tilde g)\la(\tilde g\inv m \tilde g; \tilde h).
\endaligned
$$
\item"iii)" For $m=(x, m, x), n=(x, n, x)\in \cN$, we have
$$\aligned
\tilde \la(m, n)&=\la(x; m, n)=\mu(x; n, n\inv mn)\mu^*(x; m, n)\\
&=\tilde \mu(n, n\inv mn)\overline{\tmu(m, n)}.
\endaligned
$$
\item"iv)" For each $m=(y, m, y)\in \cN, \tilde g=(y, \tilde g, x)\in \cG$,
we have
$$\aligned
\tla(m; \tilde g)\tla(n; \tilde g)\overline{\tla(mn; \tilde g)}
&=\la(y; m; \tilde g)\la(y; n, \tilde g)\overline{\la(y; mn; \tilde g)}\\
&=\mu(y\tilde g; g\inv m g; g\inv n g)\overline{\mu(y; m, n)}\\
&=\tmu(\tilde g\inv mg; \tilde g\inv n\tilde g)\overline{\tmu(m, n)}.
\endaligned
$$
\endroster
Therefore the pair $(\tla, \tmu)$ is a characteristic cocycle for $\{\cG,
\cN\}$ with values in $\T$. The map: $(\la, \mu)\in \tZ_\a(\wtH, L,
A)\mapsto (\tla, \tmu)\in \tZ(\cG, \cN, \T)$ is an injective homomorphism.

Conversely, it is clear that if $(\tla, \tmu)\in \tZ(\cG, \cN, \T)$
is given, then
the pair
$(\la, \mu)$ defined by:
$$\aligned
\la(x; m; \tilde g)&=\tla((x, m, x); (x, \tilde g, x\tilde g))\\
\mu(x; m, n)&=\tmu((x, m, x); (x, n, x))
\endaligned\quad\text{for } x\in X, \tilde g\in \wtH, m, n \in L,
$$
is an element of $\tZ_\a(\wtH, L, A)$.

Finally, each cochain $c: m\in L\mapsto c(m)\in A$ gives rise to a $\T$-valued
cochain:
$$
\tilde c(x, m, x)=c(x; m)\in \T
$$
The correspondences $c\mapsto \tilde c$ and $(\la, \mu)\mapsto (\tla, \tmu)$
intertwines the respective coboundary operations. Accordingly, $\tB_\a(\wtH, L,
A)$ corresponds exactly to $\tB(\cG, \cN, \T)$. This completes the proof.
\qed
\enddemo
\proclaim{Lemma 2.2} There exists a model
construction{\rm:}
$$
(\la, \mu)\in \tZ_\a(\wtH, L, A)\mapsto \{\sM_0(\la,\mu),
\wtH,\b^{\la, \mu}\}
$$
such that
\roster
\item"i)" $\sM_0(\la, \mu)$ is a separable von Neumann algebra of
type {\twoone}{\rm;}
\item"ii)" The restriction of the action $\b^{\la, \mu}$ of $\wtH$ to
the center of
$\sM_0(\la, \mu)$
is conjugate to the covariant system $\{\sC, \wtH, \a\}${\rm;}
\item"iii)" $L=(\b^{\la, \mu})\inv(\Int(\sM(\la, \mu)));$
\item"iv)" There exists a map $u_0: m\in L\mapsto u_0(m)\in\sU(\sM)$ such that
$$\allowdisplaybreaks\aligned
&u_0(m)u_0(n)=\mu(m, n)u_0(mn), \quad m, n \in L;\\
&\b_{\tilde g}(u_0(g\inv m g))=\la(m; \tilde g)u_0(m), \quad \tilde g=(g, s)\in
\wtH=H\times \R.
\endaligned
$$
\endroster
If the original group $H$ is amenable, then the construction yields
that the factor
$\sM_0(\la, \mu)$ is necessarily {\AFD}.
\endproclaim

\demo{Proof} First, we have the corresponding characteristic cocycle
$(\tla, \tmu)\in
\tZ_\a(\cG,\allowmathbreak \cN, \T)$. Let $\sR_0$ be an AFD factor of
type {\twoone}
and $\b$ a free action of the group $\wtH=H\times \R$ on $\sR_0$. For
each $x\in X$,
consider the
$\mu_x$-twisted crossed product
$$
\sR_{\la, \mu}(x)=\sR_0\rtimes_{\b, \mu_x}L, \quad x\in X,
$$
where $\mu_x(m, n)=\tmu(x; m, n), m, n \in L$. Let $\{u(x; m): m\in
L\}, x\in X,$ be the
$\mu_x$-projective unitary representation of $L$ in $\sR_{\la, \mu}(x)$
associated with the twisted
crossed product. Then we have
$$
\sR_{\la, \mu}(x)=\sR_0\vee \{u(x; m): m\in L\}\dprime, \quad x\in X.
$$
For each $\tilde g=(y, \tilde g, x)\in \cG$, set
$$\aligned
\b_{\scriptscriptstyle{(y, \tilde g, x)}}^{\la, \mu}(a)&=\b_{\tilde
g}(a)\quad \text{for}\quad a\in
\sR_0;\\
\b_{\scriptscriptstyle{(y, \tilde g, x)}}^{\la, \mu}(u(x;
m))&=\tla(y; gmg\inv; y, \tilde g, x)u(y;
gmg\inv), \quad m\in L.
\endaligned
$$
It is routine to check that $\b^{\la, \mu}$ is indeed an action of
the groupoid $\cG$ on the Borel
field $\{\sR_{\la, \mu}(x): x \in X\}$ of factors of type {\twoone}
and the {\vna}:
$$
\sM_0(\la, \mu)=\int_X^\oplus \sR_{\la, \mu}(x)\txd\txm(x)
$$
accommodates the required action $\b^{\la, \mu}$ of $\wtH$.

If $H$ is amenable, then $L$ is also amenable, which makes each
$\sR_{\la, \mu}(x)$
approximately finite dimensional, and therefore $\sM_0$ is AFD.
\qed
\enddemo

\proclaim{Theorem 2.3} Let $\{\sC, \R, \th\}$ be an ergodic flow and
$\a$ an action of a
countable discrete group $H$ on the flow, i.e., $\a$ is a
homomorphism of $H$ into the
group $\Aut_\th(\sC)$ of automorphisms of $\sC$ commuting with the
flow $\th$. Let
$L$ be a normal subgroup of $H$ contained in the kernel $\Ker(\a)$ of
$\a$, i.e., $L$
does not act on $\sC$ at all. Consider the product group
$\wtH=H\times \R$ whose joint
action on $\sC$ will be denoted by $\a$. Then we have a functorial model
construction{\rm:}
$$
(\la, \mu)\in \tZ_\a(\wtH, L, \sU(\sC))\mapsto \{\tM(\la, \mu),
\tilde \a^{\la, \mu}\}
$$
such that
\roster
\item"i)" The restriction of the covariant system $\{\tM(\la, \mu),
\wtH, \tilde \a^{\la, \mu}\}$
   to the center is conjugate to $\{\sC, \wtH, \a\};$
\item"ii)" The {\vna} $\tM(\la, \mu)$ is type {\twoinf} and admits a
faithful semi
finite normal trace $\tau$ which is transformed in the following fashion{\rm:}
$$
\tau\scirc \a_{g, s}=e^{-s}\tau, \quad (g, s)\in \wtH;
$$
\item"iii)" $L=(\tilde \a^{\la, \mu})\inv(\Int(\tM(\la, \mu)))$ and
$\tM(\la, \mu)$ admits
a map $u^{\la, \mu}: m\in L\mapsto u^{\la, \mu}(m)\in \sU(\tM)$ such that
$$\allowdisplaybreaks\aligned
&u^{\la, \mu}(m)u^{\la, \mu}(n)=\mu(m, n)u^{\la, \mu}(mn), \quad m, n \in L;\\
&\tilde\a_{g, s}^{\la, \mu}(u^{\la, \mu}(m))=\la(gmg\inv; g,
s)u^{\la, \mu}(gmg\inv),
\quad (g, s)\in \wtH;
\endaligned
$$
\item"iv)" To each $f\in \Map(L, \sU(\sC))$, there corresponds an isomorphism
$\sig_f: \tM(\la, \mu)\mapsto \tM(\la\part_1f, \mu\part_2f)$ such that
$$\allowdisplaybreaks\aligned
&\sig_f(u^{\la, \mu}(m))=f(m)u^{\la,\mu}(m), \quad m\in L;\\
&\sig_f\scirc \tilde\a_{g,s}^{\la, \mu}\scirc
\sig_f^{-1}=\tilde\a_{g, s}^{\la\part_1f,
\mu\part_2f},
\quad (g, s)\in \wtH;\\
&\sig_{f_1}\scirc \sig_{f_2}=\sig_{f_1f_2}, \quad f_1, f_2 \in
\Map(L, \sU(\sC));\\
&\tau\scirc \sig_f=\tau, \quad f\in \Map(L, \sU(\sC)).
\endaligned
$$
\endroster
Consequently, the restriction $\a^{\la, \mu}$ of $\tilde \a^{\la,
\mu}$ to the fixed point
subalgebra{\rm:}
$$
\sM(\la, \mu)=\tM(\la, \mu)^\th
$$
is an action of the group $H$ on a factor of type {\rm
I\!I\!I} whose flow of weights is precisely $\{\sC, \R, \th\}$ such
that its modular
characteristic invariant is given by{\rm:}
$$
\chim(\a^{\la, \mu})=[\la, \mu]\in \La_\a(\wtH, L, \sU(\sC)).
$$
If the group $H$ is amenable, then the factor $\sM$ of the model
$\{\sM(\la, \mu), H,
\a^{\la, \mu}\}$ is approximately finite dimensional.
\endproclaim
\demo{Proof} We continue the discussion in the proof of Lemma 2.2.
Let $\sR_{0,1}$
be an AFD factor of type {\twoinf} equipped with a one parameter
automorphism group
$\{\th_s: s \in \R\}$ scaling trace, i.e., $\tau\scirc
\th_s=e^{-s}\tau, s\in \R,$ with $\tau$
a faithful semi finite normal trace on $\r01$. We set
$$\aligned
\tM(\la, \mu)=\sM_0(\la, \mu)\botimes \r01&\quad\text{and}\quad
\tM(\la, \mu; x)
=\sR_{\la,\mu}(x)\botimes \r01;\\
\tM(\la, \mu)&=\int_X^\oplus \tM(\la, \mu; x)\txd\txm(x).
\endaligned
$$
Let $\d_\txm$ be the modulus of the quasi-invariant measure $\txm$ on
$X$ relative to
the groupoid
$\cG$, i.e.,
$$
\d_\txm(y, \tilde g, x)=\frac{\txd(\txm\scirc \tilde
g\inv)}{\txd\txm}\Big(x\Big), \quad (y, \tilde g, x) \in \cG,
$$
and set
$$
\rho_\txm(\tilde g)=s+\log(\d_\txm(\tilde g)), \quad \tilde g\in \cG.
$$

We then define an action $\{\as{y, \tilde g, x}{\la, \mu}\}$ of $\cG$
on the field
$\{\tM(\la, \mu; x): x \in X\}$ in the following fashion: with
$\tilde g=(y, \tilde g, x),
\tilde g=(g, s)\in \wtH=H\times \R$,
$$\CD
\tM(\la, \mu; x)=\sR_{\la, \mu}(x)\botimes \r01@>\bs{y,
\tilde g, x}{\la, \mu}\otimes \th_{\rho_\txm(y,
\tilde g, x)}>\as{y, \tilde g, x}{\la, \mu}> \tM(\la, \mu; y)=\sR_{\la,
\mu}(y)\botimes \r01.
\endCD
$$
Then we have for each $a\in \tM(\la, \mu)$ and $\tilde g=(g, s)\in \wtH$
$$\allowdisplaybreaks\aligned
\tau\scirc \a_{\tilde g}^{\la, \mu}(a)&=\int_X \tau_y\Big(\big(\a_{\tilde
g}^{\la, \mu}(a)\big)(y)\Big)\txd
\txm(y) =\int_X\tau_y\Big(\as{y, \tilde g, x}{\la,\mu}(a(x))\Big)\txd\txm(y)\\
&=\int_X\frac{\txd\tau_y\scirc \a_{(y, \tilde g,
x)}}{\txd\tau_x}\tau_x(a(x))\d_\txm(y, \tilde g, x)\txd\txm(x)\\
&=\int_Xe^{-\rho_\txm(y, \tilde g, x)}\tau_x(a(x))\d_\txm(y, \tilde
g, x)\txd\txm(x)\\
&=\int_Xe^{-s-\log(\d_\txm(y, \tilde g, x))}\tau_x(a(x))
\d_\txm(y, \tilde g, x)\txd \txm(x)\\
&=\int_Xe^{-s}\tau_x(a(x))\txd\txm(x)=e^{-s}\tau(a).
\endaligned
$$
Therefore the trace transformation rule of the integrated action
$\ta^{\la, \mu}$ of
$\wtH$ on $\tM(\la, \mu)$:
$$
\tau\scirc \ta_{g, s}=e^{-s}\tau, \quad (g, s)\in \wtH,
$$
follows. The rest of the assertions in the theorem follow easily now.
This completes the
proof.
\qed
\enddemo

\proclaim{Corollary 2.4} Suppose that $G$ is a countable discrete group with a
preassigned normal subgroup $N$  and that
$\{\sC, \R, \th\}$ is an ergodic flow. Let $\a$ be an action of $G$
on the flow $\{\sC,
\R, \th\}$ such that $N\i \Ker(\a)$, i.e., a homomorphism of $G$ into the group
$\Aut_\th(\sC)$ of automorphisms commuting with $\th$ which maps $N$ to the
identity. So the action $\a$ factors through the quotient group $Q=G/N$. Fix a
cross-section $\fs\!: Q\mapsto G$ of the quotient map $\pi\!\!: g\in
G\mapsto gN\in Q$.
   Then for any modular obstruction cocycle $[c, \nu]\in \thasout(\wtG,
N, \sU(\sC))$, there
exists an outer action $\a$ of $G$ on an infinite factor $\sM$ with
$N=\a\inv(\cntr(\sM))$ such that the associated modular obstruction
$\Obm(\a)$ is
precisely the cohomology class $\Obm(\a)=[\xi, \nu]\in \thasout(\wtG,
N, \sU(\sC))$.

If $G$ is amenable, then the construction of $\{\sM, G, \a\}$ yields
that $\sM$ is
approximately finite dimensional.
\endproclaim
\demo{Proof} Denote the unitary group $\sU(\sC)$ by $A$ for short and
$\wtH=H\times \R$. Let $[c]=\part([\xi, \z])\in \ththr(G, \T)$ be the
image of the
cohomology class $[\xi, \nu]$ under the map $\part: \thasout(\wtG, N, A)\mapsto
\ththr(G, \T)$ of
\cite{KtT2: Lemma 2.11}. Consider the resolution system:
$$\CD
1@>>>M@>>>H(c)@>\pig>>G@>>>1
\endCD
$$
of the cocycle $c\in\tzthr(G, \T)$. As $\pi_G^*([c])=1$ in $\ththr(H,
\T)$, we have
$\Inf([\xi, \nu])=\pi_G^*\scirc \part([\xi, \nu])=1\in \ththr(H, \T)$
in the modified
{\hjr}, so that we can find $\chi=[\la, \mu]\in \La_\a(\wtH, L, A)$
such that $[\xi,
\nu]=\d(\chi)$ by \cite{KtT2: Theorem 2.7}. Now Theorem 2.3 yields the
existence of a
covariant system $\{\sM, H, \a^{\la, \mu}\}$ such that $L=(\a^{\la,
\mu})\inv(\cntr(\sM))$ and $\chi(\a)=\chi$. With $\sh$ a
cross-section of $\pig$, we set
$$
\a_g=\a_{\sh(g)}^{\la, \mu} \quad \text{on}\quad \sM, \ g\in G,
$$
to obtain the outer action $\a$ of $G$ on $\sM$ such that
$\Obm(\a)=[\xi, \nu]$.
This completes the proof.
\qed
\enddemo

\head{\bf \S3. Reduction of Invariants for the Case of Type
I\!I\!I$_{\pmb\la}$,
$\pmb{0<\la<1.}$}
\endhead

First, fixing $0< \la<1$, we set
$$T=-2\pi/(\log \la)>0\quad \text{and}\quad T'=-\log \la>0.\tag3.1
$$
Let $\sC=L^\infty(\R/T'\Z)$ and $A=\sU(\sC)$. The action
$\th$ of the real line $\R$ on $\T=\R/T'\Z$ is by translation.

\proclaim{Lemma 3.1} In the above context, the first cohomology group
$\thth1(\R,
A)$ has the following structure{\rm:}
\roster
\item"i)" $\thth1(\R, A)\cong \rtz;$
\item"ii)" The following exact sequence splits{\rm:}
$$\CD
1@>>>\tbth1(\R, A)@>>>\tzth1(\R,
A)@>>\underset{\sz}\to\longleftarrow>\rtz@>>>0
\endCD
$$
\endroster
\endproclaim
\demo{Proof} i) This follows from Proposition 1.8 by setting
$Y=\{0\}\in X=\rt'z$.
Because
$$
\cG_Y=\{0\}\times T'\Z,\ A(0)=\T\ \text{ and }\ \thth1(T'\Z,
\T)=\widehat{T'\Z}\cong\rtz.
$$

ii) The notation $[x], x\in \R$, stands for the Gauss symbol,
i.e., the unique integer $n$ such that $n\leq x<n+1$.
Set 
$$
c(s, t, x)=\exp\left( \txti T's\left(\biggl[\frac{x}{T'}\biggr]-\left[\frac
{x+t}{T'}\right]\right)\right), \quad s, t, x \in \R,
\tag3.2
$$
and observe that the function $c(s, t, \cdot)$ is periodic in the
last variable $x$ with
period $T'$, so that it can be viewed as a function over $X=\rt'z$.
Furthermore, we
have
$$\allowdisplaybreaks\aligned
c(r+s, t, x)&=c(r, t, x)c(s, t, x);\\
c(r, s+t, x)&=c(r, s, x)c(r, t, x+s);\\
c(T, t, x)&=1.
\endaligned
$$
Thus $\{c(s, \cdot, \cdot)\}$ is a one parameter subgroup of
$\tzth1(\R, A)$ with
period $T$. Hence the map: $\dst\in \rtz\mapsto c(\dst,
\cdot, \cdot)\in \tzth1(\R, A)$ is a cross-section $\sz$ of the map
$\piz:\tzth1(\R,
A)\mapsto \thth1(\R, A)=\rtz.$
\qed
\enddemo

\proclaim{Corollary 3.2} If $\sM$  is a separable factor of type {\threel},
$0<\la<1$, then the association of extended modular automorphism to
each $\ds\in \rtz$ with $T=-2\pi/(\log\ \la)$ and $\f\in \fW_0(\sM)${\rm:}
$$
(\ds, \f)\in \rtz\times \fW_0(\sM)\mapsto \sig_{\sz(\dst)}^\f\in \cntr(\sM)
$$
is equivariant in the sense that if $\a$ is an isomorphism of $\sM$ 
onto another
factor $\sN$, then
$$
\a\scirc\sig_{\sz(\dst)}^{\f}\scirc \a\inv=\sig_{\a(\sz(\dst))}^{\a(\f)},
\quad \dst\in\rtz.\tag3.3
$$
In fact, for each $\f\in\fW_0(\sM)$, there exists a non-singular positive $H\in
\sC_\f$ such that with $\rho=\f(H\ \cdot\ )$,
$$
\sig_s^\rho=\sig_{\sz(\dst)}^\f, \quad s\in \R,
$$
where $\dst=s+T\Z\in \rtz$.
\endproclaim
\demo{Proof} The last lemma shows that $\{\sz(\dst)\}$ is a one parameter
subgroup of $\tzth1(\R, A)$ with period $T$. Hence for each $\f\in
\fW_0(\sM)$,
$\{u_\f(\dst): s \in \R\}=\{b_\f(\sz(\dst)): s\in \R\}$ is a one
parameter subgroup
of $\tsU(\sM)\cap \sD^\f$ with period $T$ (for notations, see
\cite{Tk2: Page 458}) and $\{\sig_{\sz(\dst)}^\f: s\in \R\}$ is a one parameter
automorphism group of $\sM$ with period $T$ which leaves $\sM_\f$
pointwise invariant.

Fix $\f\in \fW_0(\sM)$ and take a non-singular positive $k\in \sC_\f$ such that
$\Ad(k^{-iT})\allowmathbreak=\sig_T^\f$, so that $\p=\f(k\ \cdot)$ is
a {\fwt} with
$\sig_T^\p=\id$. Then it follows that $\p^{\txti T}$ belongs to the
center $\sC$, (in
fact it generates $\sC$). With $h\in \sC$ such that $h^{\txti T}=\p^{\txti T}$,
we have
$$
\partth(h^{-is}k^{is}\f^{is})=\partth(h^{-is}\p^{is})=\sz(\dst),
\quad s \in \R.
$$
Hence we get
$$
\partth(u_\f(\dst))=\sz(\dst)=\partth(h^{-is}k^{is}\f^{is}), \quad s\in \R.
\nopagebreak
$$
This means that $v(s)=u_\f(\dst)h^{is}k^{-is}\f^{-is}\in \sU(\sM), s
\in \R$. Now as $k\in
\sC_\f$, $v(s)$ and $\sD_\f$ commute, consequently so do $\f$ and $v(s)$. Hence
$\{v(s)\}$ is a periodic one parameter unitary group  in the
centralizer $\sM_\f$ with period
$T$ such that
$$\aligned
\sig_{\sz(\dst)}^\f&=\Ad(v(s))\scirc \sigps, \quad s\in \R.
\endaligned
$$
Since both $\sig_{\sz(\dst)}^\f$ and $\sigps$ leave $\sM_\f$
pointwise invariant,
$\{v(s)\}$ is contained in the center $\sC_\f$. Thus there exists a
non-singular $\ell\in
\sC_\f$ such that $v(s)=\ell^{is}$. Therefore the operator
$H=k\ell\in \sC_\f$ gives
a {\fwt} $\rho$ such that $\sigrs=\sig_{\sz(\dst)}^\f, s \in \R$, as required.

The equivariance of $\sig_{\sz(\dst)}^\f$ follows from that of the
Falcone-Takesaki cross-section $b_\f$:
$$\allowdisplaybreaks\aligned
\sig_{\sz(\ds)}^{\f}&=\tAd(b_{\f}(\sz(\ds))),\quad  \ds \in \rtz;\\
\a(b_\f&(\sz(\ds))=b_{\f\scirc\a\inv}(\a(\sz(\ds))).
\endaligned
$$
This completes the proof.
\qed
\enddemo

\proclaim{Theorem 3.3} If $\sM$ is a separable factor of type {\threel}, then
the intrinsic invariant $\Theta(\sM)$ lives in the group{\rm:}
$$
\La(\Aut(\sM)_\txm, \cnt(\sM), \T)\cong
\La_{\mod\times\th}(\Aut(\sM)\times\R,\cntr(\sM), \sU(\sC)),
$$
where
$$
\Aut(\sM)_\txm=\{(\a, s)\in \Aut(\sM)\times \R: \mod(\a)=\dstp\in\rt'z\},
$$
and  $\{\sC, \R, \th\}$ is the flow of weights on $\sM$, i.e.,
$\sC=L^\infty(\rt'z)$ and $\th$ is the translation with $T'=-\log\la.$
\endproclaim

\demo{Proof} This is  simply restating Proposition 1.11 with $Y=\{0\}\in
X=\rt'z$. We leave the detail to the reader.
\qed
\enddemo

Before going any further, we need to discuss the structure of the
third cohomology
group $\ththr(\Qm, \T)$. So let $\Qm$ be a  discrete group equipped with a
distinguished torsion free central element $z_0$. We denote by $Z$
the cyclic subgroup
generated by $z_0$, which is isomorphic to the integer group $\Z$
under the map:
$n\in \Z\mapsto z_0^n\in Z$. We denote the quotient group $\Qm/Z$ by $Q$. The
elements of $\Qm$ are denoted by $\tp, \tq, \tr$ and so on and their
quotient images
will be denoted by plain $p, q, r$ and so on. The quotient map: $\Qm
\mapsto Q$ will be
denoted by $\pim$, so $p=\pim(\tp), q=\pim(\tq), r=\pim(\tr)$.
A cross-section $\fs_\txm$ of the exact sequence:
$$\CD
1@>>>Z@>>>\Qm@>\pim>\underset{\fs_\txm}\to\longleftarrow>Q@>>>1
\endCD
$$
gives rise to an element $n\in \tztw(Q, \Z)$ such that
$$
\fs_\txm(p)\fs_\txm(q)=z_0^{n(p, q)}\fs_\txm(pq), \quad p, q \in  Q.
$$

{\smc Definition 3.4. } A cocycle $c\in \tzthr(\Qm, \T)$ is  said to
be {\it standard}
if there exists a $d_c\in \tC^2(Q, \T)$ such that
$$
c(\tp z_0^m, \tq z_0^n, \tr z_0^\ell)=d_c(q, r)^m c(\tp, \tq, \tr)
$$
for each $m, n, \ell  \in \Z$ and $\tp, \tq, \tr \in \Qm$. We denote
the group of all standard
3-cocycles by $\tZ_\txs^3(\Qm, \T)$. The cochain $d_c$ is called the {\it
$d$-part} of the standard cocycle $c\in \tZ_\txs^3(\Qm, \T)$.

Needless to say, the above definition makes sense only when we fix
the distinguished
element $z_0$ in the center of $\Qm$. This element will stand for the
automorphism
$\th_{T'}$ on the discrete core $\tM_d$ of a factor $\sM$ of type
{\threel} which scales the
trace by
$\la$.

\proclaim{Theorem 3.5} {\rm i)} Each $c\in \tzthr(\Qm, \T)$ is
cohomologous to a
standard cocycle $c_\txs\in \tZ_\txs^3(\Qm, \T)$.

{\rm ii)} The d-part $d_c$ of a standard cocycle $c\in
\tZ_\txs^3(\Qm, \T)$ is necessarily a cocycle, i.e., $d_c\in \tztw(Q, 
\T)$ and c-part $c(\tp,\tq,\tr)$
satisfies\;
$$
 c(z_0,\tq, \tr)=d_c(q,r).
$$
{\rm iii)} A standard cocycle  $c\in \tzthr(\Qm, \T)$
is a coboundary if and only if
there
is a function $f\in \tC^2(\Qm, \T)$ with $c=\part_\Qm f$ satisfying
$$\allowdisplaybreaks\left\{\aligned
f&(\tp z_0^m, \tq z_0^n)=f(z_0, q)^m
f(\tp,\ \tq);\\
 f&(z_0^n,\ q)f(z_0^m,\ q)=f(z_0^{m+n},\ q);\\
f&(\tp,\ 1)=f(1,\ \tq)=1.
\endaligned\right.
$$
for all $\tp,\ \tq\in \Qm$.

More precisely,
there exists an element $f\in
\tC^2(\Qm,\T)$ such that the function $f(\tp, \cdot)$ of the second 
variable factors through the quotient map $\pim$
and
$$\allowdisplaybreaks\aligned
&c(\tp,\tq,\tr)=(\part_\Qm f)(\tp,\tq,\tr);\\
&d_c(q, r)=f(z_0, r)f(z_0, q)f(z_0, qr)\inv;\\
&f(\tp z_0^m, \tq z_0^n)=f(z_0, q)^mf(\tp, \tq);\\
& f(z_0^n,\ q)f(z_0^m,\ q)=f(z_0^{m+n},\ q)
\endaligned
$$
for any pair $q, r\in Q$ and $\tq, \tr\in \Qm$ with $q=\pim(\tq), r=\pim(\tr)$.
\endproclaim

Choose and fix a cocycle $c\in \tzthr(\Qm, \T)$. Consider the {\vna}
$\sA=\ell^\infty(\Qm)$ of bounded functions over $\Qm$ and let $B$ be
the unitary
group $B=\sU(\sA)=\T^\Qm$, which generates an exact sequence of compact abelian
groups:
$$\CD
1@>>>\T@>i>>B@>j>\underset{\sj}\to\longleftarrow>C@>>>1.
\endCD
$$
The exact sequence splits so that the group $C$ is identified with
the subgroup of $B$
consisting of all $\T$-valued functions on $\Qm$ whose value at the
identity $1\in \Qm$
is $1$. The map $j$ is given by:
$$
(jf)(\tp)=\frac{f(\tp)}{f(1)}, \quad f\in B.
$$
The translation action of $\Qm$ does not leave the subgroup $C$
invariant, but the little
twisted action given by
$$
(\a_\tp f)(\tq)=\frac{f(\tq\tp)}{f(\tp)}, \quad f\in C,
$$
is consistent with the right translation action $\a$ of $\Qm$ on $B$
as seen below:
$$\allowdisplaybreaks\aligned
(j\scirc\a_\tp (f))(\tq)&=\frac{(\a_\tp f)(\tq)}{(\a_\tp
f)(1)}=\frac{f(\tq\tp)}{f(\tp)}\\
&=(\a_\tp\scirc j(f)(\tq), \quad f\in C.
\endaligned
$$
The action $\a_{z_0}$  will be denoted by $\th$. The fixed point
subgroup $C^\th$ will
be denoted by
$L$, which is the subgroup consisting of those functions  $f\in C$ such that
$$
\partth f=(\th f)f\inv=\text{Constant}\in \T,
$$
equivalently
$$
f(\tp z_0)=f(\tp)f(z_0), \quad \tp \in \Qm.
$$
Therefore, on $E=j\inv(L)$, the coboundary operator $\partth$ becomes the
evaluation of the value of a function $f\in E$ at $z_0$, i.e.,
$$
\partth f=f(z_0), \quad f\in E.
$$

The cocycle $i_*c\in \tzath(\Qm, B)$ is cobounded by the element $u=u_c\in
\tcatw(\Qm, B)$ defined by:
$$
u(x; \tp, \tq)=c(x, \tp, \tq), \quad x, \tp,\tq\in \Qm,
$$
as seen below:
$$\allowdisplaybreaks\aligned
1&=c(\tp, \tq, \tr)\overline{c(x\tp, \tq, \tr)}c(x, \tp\tq,
\tr)\overline{c(x, \tp, \tq\tr)}
c(x, \tp, \tq);\\
c(\tp&, \tq, \tr)=c(x\tp, \tq, \tr)\overline{c(x, \tp\tq, \tr)}c(x,
\tp, \tq\tr)\overline{c(x, \tp,
\tq)}\\
&=(\a_\tp u)(x; \tq, \tr)u(x; \tp\tq, \tr)\inv u(x; \tp, \tq\tr)u(x;
\tp, \tq)\inv\\
&=(\a_\tp (u(\tq, \tr))u(\tp\tq, \tr)\inv u(\tp, \tq\tr)u(\tp, \tq)\inv)(x)\\
&=(\part_\Qm u)(x; \tp, \tq, \tr).
\endaligned
$$
In fact, we have
$$
\tH_\a^n(\Qm, B)=\{1\}, \quad n=1, 2, \cdots.
$$
As $c(1, \tq, \tr)=1$, the element $u (\tq, \tr)\in B, \tq, \tr\in
\Qm,$ belongs to $C$. But to
distinguish $u(\tq,  \tr)\in C$ from $u(\tq, \tr)\in B$, we will
denote $u(\tq,  \tr)\in C$ by
$j_*u(\tq,  \tr)$. Since  $(\part_\Qm j_*u(\tp, \tq, \tr)))=j(c(\tp,
\tq, \tr))=1$, $j_*(u)$ is an
element of
$\tzatw(\Qm, C)$, which gives rise to an exact sequence:
$$\aligned
&\CD
1@>>>C@>>>C\rtimes_{\a,
j_*u}\Qm@>\dpi>\underset{\fs}\to\longleftarrow>\Qm@>>>1;
\endCD\\
&(a, \tp)(b, \tq)=(a\a_\tp(b)(j_*u)(\tp, \tq), \tp\tq), \quad a, b
\in C, \tp, \tq\in \Qm;\\
&\hskip 1in \dpi(a, \tp)=\tp,\quad \fs(\tp)=(1, \tp),
\endaligned
$$
so that
$$
 j_* (u(\tp, \tq))=\fn_C(\tp, \tq)=\fs(\tp)\fs(\tq)\fs(\tp\tq)\inv, \quad
\tp, \tq \in \Qm.
$$
The action $\a_\tp\in \Aut(C)$ of $\tp$ on $C$ is given as the
restriction of the inner automorphism $\Ad(\fs(\tp))\in
\Int(C\rtimes_{\a, j_*u}\Qm)$.
The coherence of the action $\a_\tp$ on
$B$ and $\Ad(\fs(\tp))$ on $C$ intertwined by $\dpi$ gives rise to an element
$(\la, \mu)\in \tZ_\a(C\rtimes_{\a, j_*u}\Qm, C, \T)$ such that
$$\aligned
\a_\tp(\sj(\fs(\tp)\inv a\fs(\tp)))&=\la(a, \tp)\sj(a), \quad (a,
\tp)\in C\rtimes_{\a,
j_*u}\Qm;\\
\mu(a, b)&=1, \quad a, b \in C.
\endaligned
$$

   The fixed point subalgebra $\sA^\th$ is identified with
$\sB=\ell^\infty(Q)$ and
therefore the fixed point subgroup $B^\th$ is identified with
$\sU(\sB)=\T^Q$, the
compact abelian group of all $\T$-valued functions over $Q$. Set
$K=B^\th/i(\T)\subset
L=C^\th$. We then have the following:
\proclaim{Lemma 3.6}
There exists a cross-section $\dfs\!: \Qm\mapsto C\rtimes_{\a,j_*(u)}\Qm$ such that the element $z=\fs(z_0)$ commutes with
$\dfs(\Qm)$. Hence the
associated cocycle $\dot\fn_L=\part_\Qm \dfs$ takes values in $L$ and therefore
$\dot\fn_L\in\tzatw(\Qm, L)$.
\endproclaim
\demo{Proof}
Since $\ththr(\Z, \T)=\{1\}$, we  may assume that $c|_{\Z^3}=1$. Hence the
restriction $u|_{\Z^2}$ is a cocycle. Hence  there exists  $v\in\tC^1(\Z,  B)$
such that $u|_{\Z^2}=\part_\Z v$. Extending $v$ to $\Qm$ and replacing
$u$ by $(\part_{\Qm}  v)\inv u$, we may and do assume that the restriction
$u|_{\Z^2}$ of $u$ to $\Z^2$ is trivial, i.e., $u(z_0^m, z_0^n)=1$. Hence
we have $\fs(z_0^m)=\fs(z_0)^m$.  Set $z=\fs(z_0)$. Now  we look at
the action $\th=\Ad(z)$ on $\fs(\tp)$:
$$\allowdisplaybreaks\aligned
\th^m(\fs(\tp))&=j(u(z_0^m, \tp))\fs(\tp z_0^m)z^{-m}\\
&=j(u(z_0^m, \tp)u(\tp z_0^m, z_0^{-m}))\fs(\tp).
\endaligned
$$
Set
$$\aligned
v(\tp, m)&=u(z_0^m, \tp)u(\tp z_0^m,  z_0^{-m}), \quad \tp\in \Qm, m\in
\Z;\\
b(\tp&, m, n)=\th_m(v(\tp, n))v(\tp, m)v(\tp, m+n)\inv.
\endaligned
$$
We claim:
$$
b(\tp, \cdot, \cdot)\in \tztw(\Z, \T).
$$
The cocycle property is obvious. So we check if $b(\tp, m, n)\in \T$:
$$\allowdisplaybreaks\aligned
j(v(\tp&, m+n))\fs(\tp)=\th^{m+n}(\fs(\tp))=\th^m(j(v(\tp, n))\fs(\tp))\\
&=j(\th^m(v(\tp, n))v(\tp, m))\fs(\tp);\\
j(\th^m(v&(\tp, n))v(\tp, m)v(\tp, m+n)\inv)=1.
\endaligned
$$
thus  $b(\tp, m, n)\in \T$. The triviality $\thtw(\Z,  \T)=\{1\}$
entails the existence
of $a(\tp, \cdot)\in \tC^1(\Z, \T)$ such that $b(\tp, \cdot, \cdot)=\part_\Z
a(\tp, \cdot)$,  i.e.,
$$\allowdisplaybreaks\aligned
\th_m(v(\tp&, n))v(\tp, m)v(\tp, m+n)\inv=a(\tp, m)a(\tp, n)a(\tp, m+n)\inv.
\endaligned
$$
Setting
$$
w(\tp, m)=a(\tp, m)\inv v(\tp, m)\in B,
$$
we get
$$
w(\tp, m+n)=\th^m(w(\tp, n))w(\tp, m)
$$
and
$$
z^m\fs(\tp)z^{-m}=j(w(\tp, m))\fs(\tp).
$$
Once  again the triviality $\tH_\th^1(\Z, B)=\{1\}$ gives the existence of
$a(\tp)\in B$ so  that $w(\tp, m)=\th^m(a(\tp)\inv)a(\tp)$. Hence we get
$$\aligned
z^m\fs(\tp)z^{-m}&=j(w(\tp, m))\fs(\tp)=j(\th^m(a(\tp)\inv)a(\tp))\fs(\tp);\\
&z^mj(a(\tp))\fs(\tp)z^{-m}=j(a(\tp))\fs(\tp).
\endaligned
$$
Therefore, $z$ and $\dfs(\tp)=j(a(\tp))\fs(\tp), \tp\in \Qm,$ commute.
\qed
\enddemo

\demo{Proof of {\rm Theorem 3.5}} i) With the last lemma, we replace
the cochain
$u\in
\tC^2(\Qm, B)$ by the cochain
$\dot u\in \tC^2(\Qm, B)$ given by:
$$
\dot u(\tp, \tq)=a(\tp)\a_\tp(a(\tq))u(\tp,
\tq)a(\tp\tq)\inv=(u\part_\Qm  a)(\tp, \tq),
\quad \tp, \tq\in \Qm,
$$
which still cobounds the cocycle $i_*(c)\in \tzthr(\Qm, B)$. Now as
$$
\dot \fn_L(\tq, \tr)=j(\dot u(\tq, \tr))\in L,
$$
$\dot u(\tq, \tr)\in E=j\inv(L)\subset B$ and
$$\left.\aligned
\dot u(xz_0; \tq, \tr)&=\dot u(x; \tq, \tr)\dot u(z_0; \tq, \tr);\\
\partth (\dot u(\tq, \tr))&=\dot u(z_0; \tq, \tr)\in \rtz;\\
\th^m(\dot u(\tq, \tr))&=\dot u(z_0; \tq, \tr)^m \dot u(\tq, \tr).
\endaligned\right\}\quad x, \tq, \tr\in \Qm.
$$
The calculation:
$$\allowdisplaybreaks\aligned
\dfnl(\tq z_0^m,  \tr z_0^n)&=\dfs(\tq z_0^m)\dfs(\tr
z_0^n)\dfs(\tq\tr z_0^{m+n})\inv\\
&=\dfs(\tq)z^m\dfs(\tr)z^n(\dfs(\tq\tr)z^{m+n})\inv\\
&=\dfs(\tq)\dfs(\tr)\dfs(\tq\tr)\inv\\
&=\dfnl(\tq, \tr),
\endaligned
$$
shows that
\roster
\item"i)" the cocycle $\dfnl$ is the pullback $\pims(\fn_L^Q)$ of a
cocycle $\fn_L^Q\in
\tzatw(Q, L)$;
\item"ii)" there exists  $f\in \tC^2(\Qm, \T)$ such that
$$
\dot u(\tp, \tq)=f(\tp, \tq)\sj(\fnlq(p, q)), \quad \tp, \tq \in \Qm,
p=\pim(\tp),  q=\pim(\tq).
$$
\endroster
   Therefore with $v(p, q)=\sj(\fnlq(p, q))$ we have
$$\allowdisplaybreaks\aligned
c(\tp, \tq, \tr)&=(\part_\Qm \dot u)(\tp, \tq, \tr)
=(\part_\Qm f)(\tp, \tq, \tr)(\part_\Qm \pims (v))(\tp, \tq, \tr)\\
&=(\part_\Qm f)(\tp, \tq, \tr)\a_\tp(v(q, r))v(p, qr)
\{v(p, q)v(pq, r)\}\inv\\
&=(\part_\Qm f)(\tp, \tq, \tr)\a_\tp(v(q, r))v(p, qr)
\{v(p, q)v(pq, r)\}\inv\\
&=(\part_\Qm f)(\tp, \tq, \tr)c_\txs(\tp, \tq, \tr)
\endaligned
$$
where
$$
c_\txs(\tp, \tq, \tr)=\a_\tp(v(q, r))v(p, qr)
\{v(p, q)v(pq, r)\}\inv=(\part_\Qm \pims(v))(\tp, \tq, \tr).
$$
Now the original cocycle $c$ is cohomologous to $c_\txs$.
We then obtain the following:
$$\allowdisplaybreaks\aligned
c_\txs(\tp &z_0^m, \tq z_0^n, \tr z_0^\ell)=\a_{\tp}\scirc \th^m(v(q,
r))v(p, qr)
\{v(p, q)v(pq, r)\}\inv\\
&=\a_\tp(\th^m(v(q, r))v(q, r)\inv)c_\txs(\tp, \tq, \tr)\\
&=v(z_0; q, r)^mc_\txs(\tp, \tq, \tr).
\endaligned
$$
Hence $d_c(q, r)=v(z_0; q, r)$ gives the $d$-part of $c_\txs$.

ii) Suppose that $c\in \tzthr(\Qm, \T)$ is standard.
 Then we have  $$d_c(q,r)=c(z_0,\tq,\tr),\quad c(z_0, \tq,\tr)=c(z_0\tp,\tq,\tr)\overline{c(\tp,\tq,\tr)}.$$
Therefore we compute for $\tp,\tq,\tr \in Q_m$:
$$\allowdisplaybreaks\aligned
d_c(q,r)&\overline{d_c(pq,r)}d_c(p,qr)\overline{d_c(p,q)}\\
=&c(z_0,\tq,\tr)\overline{c(z_0,\widetilde{pq},\tr)}c(z_0,\tp,\widetilde{qr})\overline{c(z_0,\tp,\tq)}\\
=&c(z_0,\tq,\tr)\overline{c(z_0,\tp\tq z_0^{-n(p,q)},\tr)}c(z_0,\tp,\tq\tr z_0^{-n(q,r)})\overline{c(z_0,\tp,\tq)}\\
=&c(z_0,\tq,\tr)\overline{c(z_0,\tp\tq,\tr)}c(z_0,\tp,\tq\tr )\overline{c(z_0,\tp,\tq)}\\
=&\overline{c(\tp,\tq,\tr)}c(z_0\tp,\tq,\tr)\overline{c(z_0,\tp\tq,\tr)}c(z_0,\tp,\tq\tr )\overline{c(z_0,\tp,\tq)}\\
=&\partial_{Q_m}(\bar c)(z_0,\tp,\tq,\tr)=1.\\
\endaligned
$$
Hence we conclude that $d_c\in \tztw(Q, \T)$.

iii) Suppose that $c=\part_\Qm f, f\in \tC^2(\Qm, \T)$, is standard
with $d$-part $d$.
We then compute:
$$\allowdisplaybreaks\aligned
d(q, r&)^mc(\tp, \tq, \tr)=(\part f)(\tp z_0^m, \tq z_0^n, \tr z_0^\ell)\\
&=f(\tq z_0^n, \tr z_0^\ell)\overline{f(\tp\tq z_0^{m+n}, \tr
z_0^\ell)}f(\tp z_0^m, \tq\tr
z_0^{n+\ell})\overline{f(\tp z_0^m, \tq z_0^n)}.\\
\endaligned\tag3.4
$$

Setting $\tp=1$, we get
$$\aligned
d(q, r&)^m=(\part f)(z_0^m, \tq z_0^n, \tr z_0^\ell), \quad m, n, \ell\in \Z,\\
&=f(\tq z_0^n, \tr z_0^\ell)\overline{f(\tq z_0^{m+n}, \tr
z_0^\ell)}f( z_0^m, \tq\tr
z_0^{n+\ell})\overline{f(z_0^m, \tq z_0^n)}\\
\endaligned
$$
Setting $\tq=\tr=1$, we obtain the cocycle property of $f|_{Z^2}$:
$$\aligned
1&=f(z_0^n,  z_0^\ell)\overline{f(z_0^{m+n},  z_0^\ell)}f( z_0^m,
z_0^{n+\ell})\overline{f(z_0^m, z_0^n)},\\
\endaligned
$$
Since $\tH^2(Z, \T)=\{1\}$, there exists $g\in \tC^1(Z, \T)$ such that
$$\aligned
   f(z_0^m,
z_0^n)=\overline{g(z_0^m)g(z_0^n)}g(z_0^{m+n}).
\endaligned
$$
Extend $g$ to the entire $\Qm$ and replace $f$ by $(\part_\Qm g) f$ to get
$$
f(z_0^m, z_0^n)=1, m, n \in \Z.
$$

Setting $\tq=1,\ \tr=1$ and $m=0$ in (3.4), we get
$$\aligned
1&=f(z_0^n, z_0^\ell )\overline{f(\tp z_0^{n}, z_0^\ell )}
f(\tp, z_0^{n+\ell})\overline{f(\tp, z_0^\ell)};\\
f&(\tp, z_0^{n+\ell})=((\th^n\otimes\id)f)(\tp, z_0^\ell )
f(\tp, z_0^\ell).
\endaligned
$$
Hence the cochain: $k\in \Z\mapsto f(\cdot\ , z_0^k)\in B$ is a
cocycle, thus the triviality
$\thth1(Z, B)=\{1\}$ gives the existence of $g\in \tC_\a^1(\Qm,  \T)$ such that
$$\aligned
f(\tp, z_0^\ell)=g(\tp z_0^\ell)g(\tp)\inv.
\endaligned
$$
As a constant multiplication on $g$ does not affect on the above 
identity, we may and do assume that
$g(1)=1$.
Observing
$$\allowdisplaybreaks\aligned
1&=f(z_0^n, z_0^\ell)=g(z_0^{n+\ell})g(z_0^n)\inv;\\
&g(z_0^{n+\ell})=g(z_0^n)=g(1)=1,\quad n, \ell\in Z,
\endaligned
$$
we get
$$
f(\tp, z_0^\ell)=g(z_0^\ell)\inv g(\tp z_0^\ell)g(\tp)\inv= \part_\Qm
g(\tp, z_0^\ell)\inv.
$$
Now  with $f'=f\part_\Qm g$, we obtain
$$\allowdisplaybreaks\aligned
& f'(\tp, z_0^n)=1;\quad
c=\part_\Qm f'.\\
\endaligned
$$
By replacing  $f$ by $f'$,  we may assume  $f(\tp, z_0^n)=1$ and
$c=\part_\Qm f$.

Since $f(\tp, z_0^n)=1$, for triplet $(\tp,\ z_0^n,\ \tq)$  in (3.4), we get
$$\allowdisplaybreaks\aligned
f(z_0^n,\ \tq)&\overline{f(\tp z_0^n,\ \tq)}
f(\tp,\ z_0^n\tq)\overline{f(\tp,\ z_0^n)}=c(\tp,\ z_0^n,\ \tq)=1;\\
f&(\tp z_0^n,\ \tq)=f(z_0^n,\ \tq)
f(\tp,\ z_0^n\tq)
\endaligned
$$
and for triplet $(\tp,\ \tq,\ z_0^n)$ in (3.4), we obtain
$$\allowdisplaybreaks\aligned
f(\tq,\ z_0^n)&\overline{f(\tp\tq,\ z_0^n)}f(\tp,\ \tq z_0^n)
\overline{f(\tp,\ \tq)}=c(\tp,\ \tq,\ z_0^n)=1;\\
f&(\tp,\ \tq z_0^n)=f(\tp,\ \tq)
\endaligned
$$
which implies that $f(\tp,\ \tq)$ is invariant  with
respect to $z_0^n$ in the second variable.
For triplet $(z_0^m,\ z_0^n,\ \tq)$ in (3.4), we also get
$$\allowdisplaybreaks\aligned
f(z_0^n,\ \tq)&\overline{f(z_0^mz_0^n,\ \tq)}
f(z_0^m,\ z_0^n\tq)\overline{f(z_0^m,\ z_0^n)}=c(z_0^m,\ z_0^n,\ \tq)=1;\\
f&(z_0^n,\ \tq)f(z_0^m,\ \tq)=f(z_0^{m+n},\ \tq).
\endaligned
$$
We conclude that if $c$ is a coboundary on $\Qm$, then there
is a function $f\in \tC^2(\Qm, \T)$ with $c=\part_\Qm f$ such that
for all $\tp,\ \tq\in \Qm$,
$$\allowdisplaybreaks\left\{\aligned
f&(\tp z_0^n,\ \tq z_0^m)=f(z_0^n,\ \tq)
f(\tp,\ \tq);\\
f&(z_0^n,\ \tq)f(z_0^m,\ \tq)=f(z_0^{m+n},\ \tq);\\
f&(\tp,\  1)=f(1,\ \tq)=1.
\endaligned\right.
$$
This completes the proof of the theorem.
\qed
\enddemo

We apply Theorem  3.5 to the situation:
\roster
\item"i)" the group $Q$ is the quotient group
$G/N$ of a {\dg} $G$ by a central normal subgroup $N$;
\item"ii)"  $\txm\!: G\mapsto
\rt'z$ is a homomorphism which factors through $Q$, i.e.,
$\Ker(\txm)\supset N$ so
that $\txm=\txm\scirc \pim$ employing the same notation $\txm$  for the
homomorphism from $Q$ to $\rt'z$ induced from $\txm\in\Hom(G,\allowmathbreak
\rt'z)$;
\item"iii)" the group $\Qm$ is given by:
$$
\Qm=\{\tp=(p, s)\in Q\times \R: \txm(p)=\dstp=s+T'\Z\in \rt'z\},
$$
hence the  quotient map $\pim$ is precisely the projection map
$$
\pim=\pr_1: \tp=(p, s)\in \Qm \mapsto p\in Q
$$
to the first component;
\item"iv)" the distinguished  central element $z_0$ is given by:
$$
z_0=(1, T')\in \Qm.
$$
\endroster
In this setting, we make the following:

{\smc Definition 3.7.} The {\it standard coboundaries}  $\tbsthr(\Qm, \T)$ in
$\tzsthr(\Qm, \T)$ are given by:
$$
\tbsthr(\Qm, \T)=\part_\Qm (\pims(\tctw(Q,  \T)))\subset \tbthr(\Qm,
\T).\tag3.5
$$
The {\it standard third cohomology group}  $\thsthr(\Qm, \T)$ is
defined to be the
quotient group:
$$
\thsthr(\Qm, \T)=\tzsthr(\Qm, \T)/\tbsthr(\Qm, \T),\tag3.6
$$
which is a compact abelian group.

The reason for this smaller coboundary group comes from the fact 
that when we perturb an outer
action $\a_g, g\in G$, on a factor $\sM$, we do not allow a 
perturbation by $\cntr(\sM)$ but by $\Int(\sM)$.
So even if we consider an outer action of the bigger group $\Gm$ on the 
discrete core $\tM_d$, we can not
use all of $\sU(\tM)$ but only $\sU(\sM)$.

The $d$-part $d_c$ of each $c\in\tzsthr(\Qm, \T)$ is an element of
$\tztw(Q, \T)$
and also each $\nu\in\Hom(N, \rtz)$ gives rise to an element of $\tztw(Q, \T)$:
$(q, r)\in Q^2\mapsto \nu(\fnl(p, q))\in\rtz$ under the
identification of $\dst\in \rtz$
and $e^{iT's}\in \T, s\in \R$. Hence we define $\tzsthr(\Qm, \T)*_\fs
 \Hom_G(N, \rtz)$ to
be the subgroup of $\tzsthr(\Qm, \T)\times \Hom_G(N, \rtz)$ consisting
of all those
elements $(c, \nu)$ such that
$$
d_c(q, r)=\nu(\fnl(q, r)), \quad q, r \in Q.\tag3.7
$$
Here the fiber product depends on the cocycle $\fnl\in \tztw( Q, L )$
explicitly and
therefore on the cross-section $\fs\!: Q\mapsto G$ of $\pi$.

\proclaim{Theorem 3.8} Fix $0<\la<1$ and set
$$
T=-\frac{2\pi}{\log \la},\quad T'=\frac{2\pi}T=-\log\la \quad 
\text{and}\quad X=\rt'z.
$$
Let $G$ be a group equipped with a central subgroup $N$. For a homomorphism
$\txm\!: G\mapsto X$ such that $\Ker(\txm)\supset N$, consider the  joint
action of $\wtG=G\times \R$ on $X$ given by{\rm:}
$$
T_{g, s}(x)=x - \dstp + \txm(g)\in X, \quad (g, s)\in \wtG,
$$
and the action $\a$ of $\wtG$ on $\sC=L^\infty(X)$ given by{\rm:}
$$
\a_{g, s}f(x)=f(T_{ g, s}^{-1}x), \quad x \in X=\rt'z, f \in \sC.
$$
The action of $\R$ alone denoted by $\th$ is the transitive flow with
period $T'$
and gives $\thth1(\R, A)\cong \rtz$ where $A=\sU(L^\infty(X))$.  Define $Q=G/N$
and $Q_\txm=\{(p, s)\in Q\times \R: \txm(p)=\dst\in \rt'z\}$. Define the
subgroup $\tB_{\txm, \fs}^\out(G, N, \T)$ to be the subgroup
of $\tzsthr(\Qm, \T)*_\fs \Hom_G (N, \rtz)$ consisting of all those elements
$(c, \nu)$ of the form{\rm:}
$$
c=\part_\Qm (\pims f)\quad \text{and}\quad  \nu\equiv 1 ,\tag3.8
$$
for some $f\in \tctw(Q, \T)$
and form the quotient group{\rm:}
$$
\tH_{\txm, \fs}^\out(G, N, \T)
=(\tzsthr(\Qm, \T)*_\fs \Hom_{ G}(N, \rtz))/\tB_{\txm, \fs}^\out(G, N, \T).\tag3.9
$$
   Then there is a
natural isomorphism{\rm:}
$$
\thasout(\wtG, N, A)\cong\tH_{\txm, \fs}^\out(G, N, \T).
$$
\endproclaim

The joint action of $\wtQ=Q\times \R$ on $X=\rt'z$ is transitive. But
the coboundary
group $\tB_{\a, \txs}^3(\wtQ, A)=\part_{\widetilde Q}(\tC^2(Q, A))$ is smaller than
the regular
coboundary group $\part_\wtQ(\tC^2(\wtQ, A))$ for the third cohomology group
$\ththra(\wtQ, A)$. So it is not clear whether  the straightforward
Shapiro machine
works or not. We have seen that for the relative cohomology group
$\La_\a(\wtH, L,
M, A)$ the Shapiro machinery works in Proposition 1.11. So we begin
by looking at a
resolution system $\pig: H\mapsto G$ with $L=\pi_G^{-1}(N)$ and $M=\Ker(\pig)$,
so that $c\in \tzasout(G, N, A)$ is of the form $c=\d(\la, \mu)$ for
some $(\la, \mu)\in
\tZ_\a(\wtH, L, M, A)$. Let $E=E(\la, \mu)\in \X(\wtH, L, M, A)$ be the
corresponding crossed extension equipped with a cross-section $\sj:
L\mapsto E$ such
that
$$\allowdisplaybreaks\aligned
\mu(m, n)&=\sj(m)\sj(n)\sj(mn)\inv, \quad m, n \in L;\\
\la(m, h, s)&=\a_{h, s}(\sj(h\inv mh))\sj(m)\inv,\quad (h, s)\in \wtQ.
\endaligned
$$
Let $\rho$ be the groupoid homomorphism of $\wcH=X\rtimes \wtH$ to
the stabilizer
subgroup $H_\txm=\{(h, s)\in \wtH: \txm(h)=\dstp\}$, where the map $\txm:
H\mapsto \rt'z$ is defined to be the pullback $\txm=\txm\scirc\pig$ of the
corresponding map $\txm$ of $G$ by $\pig$. The map $\rho$ is
explicitly written in
the form:
$$
\rho(\dy, h, s)=(h, s-\bracett'{\dy -\txm(h)+\ds}+\bracett'{\dy})\in
H_\txm, \quad (\dy,
h, s)\in \cH.
\tag3.10
$$
In fact, we have
$$\allowdisplaybreaks\aligned
\rho(\dy&, hk, s+t)=(hk, s+t-\bracett'{\dy -\txm(hk)+\ds+\dt}+\bracett'{\dy})\\
&=(h, s-\bracett'{\dy-\txm(h)+\ds}+\bracett'{\dy})\\
&\hskip.2in\times (k, t-\bracett'{\dy
-\txm(hk)+\ds+\dt}+\bracett'{\dy-\txm(h)+\ds})\\
&=\rho(\dy, h, s)\rho(\dy\tilde h, k, t)\quad \text{with }\tilde h=(h, s).
\endaligned
$$
\proclaim{Lemma 3.9} If $\rho$ is a homomorphism of $\wcH$ to a group $K$, then
for any $c\in \tZ^n(K, \T)$, the pullback  $\rho^*(c)$ defined by{\rm:}
$$
\rho^*(c)(\dot x; \tilde h_1, \cdots, \tilde h_n)=c(\rho(\dot x, \tilde h_1),
\rho(\dot x\tilde h_1, \tilde h_2), \cdots, \rho(\dx \tilde h_1\cdots
\tilde h_{n-1},
\tilde h_n))\tag3.11
$$
is an element of $\tZ_\a^n(\wcH, \T)$ and if  $f\in \tC^{n-1}(K, \T)$, then
$$
\part_\wtH \rho^*(f)=\rho^*(\part_K f)\tag3.12
$$
where $\rho^*(f)\in \tC^n(\wcH, \T)$ is given by
$$\aligned
(\rho^*f)&(\dx; \tilde h_1, \cdots, \tilde h_{n-1})\\
&=f(\rho(\dx, \tilde h_1), \rho(\dx\tilde h_1,
\tilde h_2), \cdots, \rho(\dx\tilde h_1\cdots \tilde h_{n-2}, \tilde h_{n-1})).
\endaligned\tag3.13
$$
Hence $\rho^*$ gives a homomorphism of $\tH^n(K, \T)$ to $\tH_\a^n(\wcH, \T)$.
\endproclaim
This follows from a direct calculation. We leave it to the reader.

\demo{Proof of Theorem {\rm 3.8}}  First, since $\rho(\dy, \tilh)=\tilh$ for
every $\tilde h\in H_\txm$, we have
$i^*\scirc
\rho^*=\id|_{\tZ^3(H_\txm, \T)}$  where $i$ is the embedding map $i_{\{0\}}$ in 
Proposition 1.8.
Next, choose $(c, \nu)\in \tzaout(\wtQ, A)$ and assume that the system $\pig:
H\mapsto G$ gives a resolution of $\part(c, \nu)\in \tzthr(G, \T)$ by
$(\la, \mu)\in
\tZ_\a(\wtH, L, M, A)$ so that $([c], \nu)=\d([\la, \mu])$ by the
modified HJR-map:
$$
\d\!: \La_\a(\wtH, L, M, A)\mapsto \thaout(G, N, A).
$$
This means that if $\sj$ is the cross-section of the crossed extension:
$$
E=E(\la, \mu)\in \X_\a(\wtH, L, M, A)
$$
associated with $(\la, \mu)$, the cocycle $(c, \nu)$ is
given by:
$$\allowdisplaybreaks\aligned
c(\tp, \tq, \tr)&=\a_\tp(\sj(\fnl(q, r)))\sj(\fnl(p,
qr))\{\sj(\fnl (p, q) \sj(\fnl(pq,
r))\}\inv;\\
\nu(n)&=[\la(\sh(n); (1, \cdot))]\in \thth1(\R, A)=\rtz,
\endaligned
$$
for $\tp=(p, s), \tq=(q, r), \tr=(r, u)\in \wtQ$ and $n\in N$ where
$\sh$ is a cross-section of the
quotient map $\pig$.  By Proposition 1.11, with $(\la_0, \mu_0)=i(\la, \mu)\in
\tZ(Q_\txm, L, M,
\T)$ we have
$$
(\la, \mu)\equiv \rho^*(\la_0, \mu_0)\quad \mod\ \tB(\wtQ, L, M, A).
$$
Therefore we may replace $(\la, \mu)$ by $\rho^*(\la_0, \mu_0)$,
i.e., we may assume that
$$\allowdisplaybreaks\aligned
\mu(\dy; m, n)&=\mu_0(m, n)\in \T, \quad m, n \in L;\\
\la(\dy; m, \tilh)&=\la_0(m, \rho(\dy, \tilh)), \quad m\in L, \
\tilh\in\wtH, \dy\in X=\rt'z.
\endaligned
$$
With $E_0=E(\la_0, \mu_0)\in \X(H_\txm, L, M, \T)$, we have
$E=\rho^*(E_0)$, i.e.,
$$
E=A\times E_0/\{(a, \bar a): a\in \T\}
$$
admits a cross-section $\sj\!: m\in L\mapsto [\fs_0(m)]\in E$ such that
$$\allowdisplaybreaks\aligned
\sj(m)\sj(n)&=\mu(m, n)\sj(mn);\\
\a_\tilh(\sj(h\inv m h))(\dy)&=\la(\dy; m, \tilh)\sj(m) (\dy) =\a_{\rho(\dy,
\tilh)}^0(\fs_0(h\inv m h)).
\endaligned
$$
Now let us compute, based on the fact that the $L$-valued cocycle
$\fnl$ does not
depend on the
$\R$-variables:
$$\allowdisplaybreaks\align
(\rho^*c^0)&(\tp, \tq, \tr; \dy)=c^0(\rho(\dy,  \tp), \rho(\dy\tp,
\tq),  \rho(\dy\tp\tq, \tr))\\
&=\a_{\rho(\dy, \tp)}^0(\fs_0(\fnl(\rho(\dy\tp, \tq), \rho(\dy\tp\tq,
\tr)))\fs_0(\fnl(\rho(\dy, \tp),
\rho(\dy\tp, \tq\tr))\\
&\hskip.2in\times
\{\fs_0(\fnl(\rho(\dy, \tp), \rho(\dy\tp,
\tq)))\fs_0(\fnl(\rho(\dy,\tp\tq),\rho(\dy\tp\tq, \tr)))\}\inv\\
&=\Big(\a_\tp(\sj(\fnl(q, r)))\sj(\fnl(p, qr))\{\sj(\fnl(p,
q))\sj(\fnl(pq, r))\}\inv\Big)(\dy)\\
&=c(\tp, \tq, \tr; \dy).
\endalign
$$
Suppose 
that $c^0\in\tbsthr(\Qm, \T)$, i.e., for some $f\in
\tctw(Q, \T)$ we have
$$\aligned
c^0(\tp, \tq, \tr)&=f(q, r)\overline{f(pq, r)}f(p, qr)\overline{f(p, q)};\\
\endaligned
$$
Then we have for $c=\rho^*c^0$:
$$\aligned
(\rho^*c^0)&(\tp, \tq, \tr; \dy)=f(\pim(\rho(\dy\tp, \tq)),
\pim(\rho(\dy\tp\tq, \tr)))\\
&\hskip .4in\times
\overline{f(\pim(\rho(\dy, \tp))\pim(\rho(\dy, \tp)),
\pim(\rho(\dy\tp, \tq\tr)))}\\
&\hskip.6in\times
f(\pim(\rho(\dy, \tp)), \pim(\rho(\dy\tp, \tq)\rho(\dy\tp\tq, \tr)))\\
&\hskip.8in\times
\overline{f(\pim(\rho(\dy, \tp)), \pim(\rho(\dy\tp, \tq)))}\\
&=f(q, r)\overline {f(pq, r)}f(p, qr)\overline{f(p, q)}\\
&=\part_\wtQ (\pr_1^* f)(\tp, \tq, \tr; \dy),
\endaligned
$$
where $\pr_1$ is  the projection map of $\wtQ$ to the
first component
$Q$,  so that $\rho^*c^0\in \part_\wtQ(\tC^2(Q, \T))\subset \part_\wtQ(\tC^2(Q, \T))$.
If
$\rho^*c^0\in \part_\wtQ(\tC^2(Q, \T))$, i.e.,
if there exists $f\in \tC^2(Q, \T)$ such that $\rho^*c^0=\part_\wtQ f$, then
for each $\tp, \tq, \tr\in Q_\txm$ we have
$$\aligned
c(\tp, \tq, \tr; \dy)&=f(q, r;\dy+\txm(p)-\ds)f(p, qr;
\dy)\overline{f(p, q; \dy)f(pq,
r; \dy)}\\
&=f(q, r; \dy)f(p, qr; \dy)\overline{f(p, q; \dy)f(pq, r; \dy)}.
\endaligned
$$
Hence we get
$$
c^0(\tp, \tq, \tr)=c(\tp, \tq, \tr; 0)=(\part_{Q_\txm}\pims
(f_0))(\tp, \tq, \tr)
$$
where $f_0(p, q)=f(p, q; 0)$. Thus we conclude $c_0\in
\part_{Q_\txm}(\pims(\tC^2(Q, \T)))$. Consequently, we get
$$
\thasth(\wtQ, A)\cong \tH_\txs^3(Q_\txm, \T).
$$

We want to compare the $d$-part $d_c$ of $c$ and the $d$-part
$d_c^0=d_{c^0}$ of $c^0$.

In  terms of $c$ and $c^0$, $d_c$ and $d_c^0$ are given by the 
following:
$$\allowdisplaybreaks\aligned
d_c(s; q, r)&=c((1, s), (q, 0), (r, 0))\\
&=\th_s(\sj(\fnl(q, r))\sj(\fnl(1, qr))\{\sj(\fnl(1, q))\sj(\fnl(q, r))\}\inv\\
&=\th_s(\sj(\fnl(q, r))\sj(\fnl(q, r))\inv\\
&=\la(\fnl(q, r), s);\\
d_c(s; q, r; \dy)&=\la(\dy; \fnl(q, r), s)=\la_0(\fnl(q, r), \rho(\dy, 1, s))\\
&=\la_0(\fnl(q, r),  s-\bracett'{\dy+\ds}+\bracett'{\dy})\\
&=\left\langle \nu(\fnn(q, r)), \Big[\frac
{y+s}{T'}\Big]-\Big[\frac  y{T'}\Big]\right\rangle \\
&=\exp\Big(\txti T'\bracettt{\nu(\fnn(q, r))}\Big(\Big[\frac 
y{T'}\Big]-\Big[\frac
{y+ s}{T'}\Big]\Big)\Big);\\
d_c^0(q, r)&= c((1, T'), (q, 0), (r, 0))\\
&=\th_{T'}(\sj(\fnl(q, r)))\sj(\fnl(q, r))\inv=\la_0(\fnl(q, r), z_0)\\
&=\dot\nu(\fnl(q, r))=\exp(-\txti T'\bracettt{\nu(\fnl(q, r))})\\ 
&=d_c(T'; q, r; 0)
\endaligned
$$
where the duality pairing of $\rtz$ and  its dual
$\Z$ will be denoted by $$\langle \ds, m\rangle=e^{-\txti T'm\bracettt{\ds}}
\ \text{for}\ m\in\Z$$ and $\ds\in \rtz$ and we write  $\dnu(m), m\in N,$ for
$e^{-\txti T'\bracettt{\nu(m)}}$ for short.

Hence $(c, \nu)$ is in $\tzasth(\wtQ, A)*_\fs\Hom_{ G}(N, \rtz)$ if and
only if $(c^0, \nu)$ is in
$\tzsthr(\Qm, \T)*_\fs \Hom_{ G}(N, \rtz)$.

Now we suppose $(c,  \nu)\in \tbasout(G, N, A)$
i.e., there exists $f\in \tcatw(Q, A)$ such that for each $\tp=(p, s),
\tq=(q, t), \tr=(r, u)\in \wtQ,$ we have
$$\allowdisplaybreaks\aligned
c(\tp, \tq, \tr; \dy)&=f(q, r; \dy\tp)f(p, qr; \dy)\overline{f(p, q; \dy)f(pq,
r; \dy)}\inv; \\
&d_c(s, q, r; \dy)=\overline{f(q, r; \dy)}f(q, r; \dy-\ds).
\endaligned
$$
Then we have, for each $\tp, \tq, \tr\in \Qm$,
$$\allowdisplaybreaks\aligned
c^0(\tp, \tq, \tr)&=c(\tp, \tq, \tr; 0)\\
&=f(q, r; 0)f(p, qr; 0)\overline{f(p, q; 0)f(pq, r; 0)};\\
d_c^0(q, r)&=\overline{f(q, r; 0)}f(q, r; 0)=1.
\endaligned
$$
Thus we get $(c^0, \nu)\in \tbmsout(G, N, \T)$. Conversely, suppose
$(c^0, \nu)\in \tbmsout(G,\allowmathbreak N, \T)$, i.e., 
 $ \nu\equiv 1 $ and for some $f\in
\tctw(Q, \T)$, 
$$\allowdisplaybreaks\aligned
c^0(\tp, \tq, \tr)&=f(q, r)\overline{f(pq, r)}f(p, qr)\overline{f(p, q)}, \quad
\quad  \tp,
\tq, \tr\in \Qm.
\endaligned
$$
Since we have seen already that $c=\rho^*c^0$ is
cobounded by $\pr_1^*f$, we have $(c,  \nu)\in \tbasout(G, N, A)$.
This completes the proof.
\qed
\enddemo

\subhead\nofrills{\bf The Map $\pmb \part$:}
\endsubhead\quad
We now want to identify the map
$$
\CD\part\!: \thasout(G, N, A)\cong\thmsout(G, N, \T)@>>> \tzthr(G,
\T)\endCD
$$
of \cite{KtT2: Theorem 2.7} in terms of $\thmsout(G, N, \T)$. To this
end, we need notations to shorten mathematical expressions. We use
$\fnz(p, q)$ for the $\Z$-valued two cocycle $\etat'(\txm(p), \txm(q)), p,
q\in Q$. We also use $\fnz(g, h)$ for $\fnz(\pi(g), \pi(h))$ omitting the
map $\pi$. Also the element $z_0=(1, T')$ appears both in $\Gm$ and
$\Qm$. So we have
$$
\bp\bq=z_0^{\fnz(p, q)}\overline{pq}; \quad \bg\bh=z_0^{\fnz(g,
h)}\overline{gh},
$$
where $\bg=(g, \bracett'{\txm(g)})\in \Gm, g\in G$ and $\bp=(p,
\bracett'{\txm(p)})\in \Qm$.

\proclaim{Lemma 3.10} Fix $(c, \nu)\in \tzsthr(\Qm,
\T)*_\fs\Hom_{ G}(N, \rtz)$. With
$$
 \txnn(g)=\fs(\pi(g))g\inv \in N, \quad  g\in G,\tag3.14
$$
and
$$
\dnu(m)=e^{\txti T'\bracettt{\nu(m)}}\in \T, \quad m\in L,
$$
set
$$
c_G(g, h, k)=\dnu( \txnn(k))^{-\fnz(g, h)}  c(\bp, \bq, \br),
\quad g, h, k\in G,\tag3.15
$$
where $p=\pi(g), q=\pi(h), r=\pi(k)$.
Then $c_G\in \tzthr(G, \T)$. The map{\rm:}
$$\CD
([c], \nu)\in \thmsout(G, N, \T)@>>>
[c_G]\in\ththr(G, \T)
\endCD
$$
is precisely the map $\part$ of \cite{KtT2: Theorem 2.7}.
\endproclaim
\demo{Proof} Since $ \txnn(g)g=\dfs(\pi(g)) , g \in G,$ where
$ \fs: Q \mapsto G $, we have for each pair
$g, h \in G$:
$$\aligned
\txnn(\pi(g), \pi(h))
=&\fs(\pi(g))\fs(\pi(h))\fs(\pi(gh))\inv\\
=&\txnn(g)g\txnn(h)h\{\txnn(gh)gh\}\inv\\
=&\txnn(g)g\txnn(h)g\inv\txnn(gh)\inv.\\
\endaligned
$$

We compute for $g, h, k, \ell \in G$ with $p=\pi(g),
q=\pi(h), r=\pi(k)$ and $s=\pi(\ell)$
$$\allowdisplaybreaks\aligned
c(\bq &, \br, \bar s)\overline{c(\overline{pq}, \br, \bar s)}c(\bp,
\overline{qr},
\bar s)\overline{c(\bp,\bq, \overline{rs})}c(\bp, \bq, \br)\\
&=c(\bq, \br, \bar s)\overline{c(z_0^{-\fnz(p,
q)}\bp\bq, \br, \bar s)}c(\bp, z_0^{-\fnz(q, r)} \bq\br, \bar s)\\
&\hskip1in\times
\overline{c(\bp, \bq, z_0^{-\fnz(r, s)}\br\bar s)}c(\bp, \bq, \br)\\
&=d_c(\pi(k), \pi(\ell))^{\fnz(g, h)}=\langle \nu(\fnn(\pi(k),
\pi(\ell))), \fnz(g, h)\rangle\\
&=\dnu\Big(\txnn(k)k\txnn(\ell)k\inv\txnn(k\ell)\inv\Big)^{\fnz(g, h)}
\endaligned
$$
and proceed to obtain the calculation:
$$\allowdisplaybreaks\align
c_G(h&, k, \ell)\overline{c_G(gh, k, \ell)}c_G(g, hk, \ell)\overline{c_G(g, h,
k\ell)}c_G(g, h, k)\\
&=\dnu(\txnn(\ell))^{-\fnz(h, k)} c(\bq,
\br, \bar s)\overline{\dnu(\txnn(\ell))^{-\fnz(gh,
k)} c(\overline{pq}, \br, \bar s)}\\
&\hskip .5in\times
  \dnu(\txnn(\ell))^{-\fnz(g, hk)} c(\bp, \overline{qr},
\bar s)\\
&\hskip1in\times
\overline{ \dnu(\txnn(k\ell))^{-\fnz(g, h)} c(\bp, \bq, \overline{rs})}\\
&\hskip1.5in\times
  \dnu(\txnn(k))^{-\fnz(g, h)} c(\bp, \bq, \bar s)\\
&= \dnu(\txnn(k)k\txnn(\ell)k\inv\txnn(k\ell)\inv)^{\fnz(g,
h)}\\
&\hskip .5in\times
  \dnu(\txnn(\ell))^{-\fnz(h, k)} \dnu(\txnn(\ell))^{\fnz(gh, k)}\\
&\hskip1in\times
  \dnu(\txnn(\ell))^ {-\fnz(g, hk)}
  \dnu(\txnn(k\ell))^{\fnz(g, h)}\\
&\hskip1.5in\times
  \dnu(\txnn(k))^{-\fnz(g, h)}\\
&=\dnu(\txnn(\ell))^{-\fnz(h, k)+\fnz(gh,k)-\fnz(g,hk)+\fnz(g,h)}\\
&=1.
\endalign
$$
Hence $c_G$ belongs to $\tzthr(G, \T)$.

Since
$$\aligned
\z_\nu(t, n)(\dot x)
&=c(\nu(n), t, \dot x)=\dnu(n)^{\left(\left[\frac{x+t}{T'}\right]-\left[\frac
{x}{T'}\right]\right)}\\
&=\left\langle \dnu(n), \left[\frac  {x+t}{T'}\right]-\biggl[\frac
{x}{T'}\biggr]\right\rangle,
\endaligned \quad  t, x \in
\R,\quad n\in N,
$$
we compute the element $a$ appeared in  (2.23) in [KtT2]  as follows;
$$\allowdisplaybreaks\align
(\th_{t}&(a(g, h))a(g, h)^*)(\dx)\\
&=\z_\nu(t; \fnn(\pi(g), \pi(h)))(\dx)\\
&\hskip.5in\times
\left\{\z_\nu(t; \txn_N(g))(\dx)
\a_g(\z_\nu(t; \txn_N(h)))(\dx)\z_\nu(t; \txn_N(gh))^*(\dx)\right\}^*\\
&=\left\langle \nu(\fnn(\pi(g), \pi(h))), 
\left[\frac {x+t}{T'}\right]-\biggl[\frac
{x}{T'}\biggr]\right\rangle\\
&\hskip1cm\times\left\{
\left\langle\nu(\txnn(g)),\left[\frac {x+t}{T'}\right]-\biggl[\frac
{x}{T'}\biggr]\right\rangle\right.\\
&\hskip1cm \times
\left\langle\nu(\txnn(h)), {\left[\frac
{x-\bracett'{\txm(g)}+t}{T'}\right]-\left[\frac
{x-\bracett'{\txm(g)}}{T'}\right]}\right\rangle\\
&\hskip1cm\times
\left.\left\langle\nu(\txnn(gh)), \biggl[\frac
{x}{T'}\biggr]-\left[\frac {x+t}{T'}\right]\right\rangle\right\}^* \\
\text{since}\ \fnn&(\pi(g), \pi(h))=\txnn(g)g\txnn(h)g^{-1}\txnn(gh)^{-1}\ 
\text{and}\ \nu \ \text{is}\ G\text{-invariant},\\
&=\dnu(\txnn(h))^{\left(\left[\left.\left.\frac
{x-\bracett'{\txm(g)}}{T'}\right]-\right[\frac 
{x+t-\bracett'{\txm(g)}}{T'}\right]\right)}
\dnu(\txnn(h))^{
\left(\left[\left.\left.\frac {x+t}{T'}\right]-\right[\frac
{x}{T'}\right]\right)}\\
&=\dnu(\txnn(h))^{\left(\left[\frac {x+t}{T'}\right]-\left[\frac
{x-\bracett'{\txm(g)}+t}{T'}\right]
\right)}
\dnu(\txnn(h))^{\left(\left[\frac{x-\bracett'{\txm(g)}}{T'}\right]-\left[\frac
{x}{T'}\right]\right)},
\endalign
$$
which shows that the cochain $a$ ought to be of the following form:
$$
a(g, h)(x)=\dnu(\txnn(h))^{
\left.\left.\left(\left[\frac
{x}{T'}\right]-\right[\frac{x-\bracett'{\txm(g)}}{T'}\right]\right)}.
$$
Now we examine the proof of \cite{KtT2: Theorem 2.7}, in particular the
proof of Lemma 2.11. The split property of the exact sequence:
$$\CD
1@>>>\tbth1@>>>\tzth1@>>\underset{\sz}\to\longleftarrow>\thth1=\rtz@>>>1
\endCD
$$
allows us to choose $f=1$ in \cite{KtT2: (2.20)}.

We  compute the map $\partial$ as follows:
$$\allowdisplaybreaks\align
(\partial c)&(g,h,k)\\
&=c(\pi(g),\pi(h),\pi(k))(x)(\partial_Ga^*)(g,h,k)(x)\\
&=c(\pi(g),\pi(h),\pi(k))(x)
\dnu(\txnn(k))^{
\left(\left[\frac{x-\bracett'{\txm(g)}-\bracett'{\txm(h)}}{T'}\right]-\left[\frac
{x-\bracett'{\txm(g)}}{T'}\right]\right)}\\
&\hskip1cm\times
\dnu(\txnn(hk))^{
\left(\left.\left.\left[\frac{x-\bracett'{\txm(g)}}{T'}\right]-\right[\frac
{x}{T'}\right]\right)}
\dnu(\txnn(h))^{
\left(\left[\frac {x}{T'}\left]
-\left[\frac{x-\bracett'{\txm(g)}}{T'}\right]\right.\right.\right)}\\
&\hskip2cm\times
\dnu(\txnn(k))^{\left(\left[\frac{x}{T'}\right]
-\left[\frac{x-\bracett'{\txm(gh)}}{T'}\right]\right)}\\
&=c(\pi(g),\pi(h),\pi(k))(x)
\exp(\txti \nu(\txnn(k))
\fn_Z(g,h))\\
&\hskip1cm\times
\dnu(\fn_N(h,k))^{\left(
\left[\frac
x{T'}\right]-\left[\frac{x-\bracett'{\txm(g)}}{T'}\right]\right)}.
\endalign
$$
Therefore we obtain the image of $c$ under the  map $\partial$ by 
evaluating at $0$:
$$\allowdisplaybreaks
\aligned
c(\bar p,\bar q,\bar r)&=c(\pi(g),\pi(h),\pi(k))(0)\dnu(\fn_N(h,k))^{-
\left[\frac{-\bracett'{\txm(g)}}{T'}\right]};\\
(\partial c)(g,h,k)&=c(\pi(g),\pi(h),\pi(k))(0)
\dnu(\txnn(k))^{-\fn_Z(g,h)}\\
&\hskip4cm\times
\dnu(\fn_N(h,k))^{-\left[\frac{-\bracett'{\txm(g)}}{T'}\right]}\\
&=c(\bar p,\bar q,\bar r)\dnu(\txnn(k))^{-\fn_Z(g,h)}.
\endaligned
$$
\qed
\enddemo

Summarizing the above arguments, we describe the modified {\hjr} of
\cite{KtT2: Theorem 2.7} in terms of cohomology groups with the
coefficient group $\T$  in the following:
\proclaim{Theorem 3.11}
There is a commutative diagram between the modified \linebreak
Huebschmann-Jones-Ratcliffe exact sequences of $\widetilde H$ and $\Qm$
{\rm:}
$$\eightpoint\CD
@.@. \ththr(G, \T)@>>\inf>\ththr(H, \T)\\
@.@.@A\partial AA@|\\
\thtw(H, \T)@>\Res >> \La(\wtH, L, M, A)@>\d >>
\thasout(G\times\R, N, A)@>\Inf>>\ththr(H, \T)\\
@| @V i^* V\Big\uparrow\rho^*V@V i^* V\Big\uparrow\rho^*V @|\\
\thtw(H, \T)@>\Res_\Qm>> \La(\Hm, L, M, \T)@>\d_\Qm>>
\thasout(G, N, \T)@>\Inf_\Qm>> \ththr(H, \T)\\
@.@.@V\partial_\Qm VV@|\\
@.@. \ththr(G, \T)@>>\inf>\ththr(H, \T)
\endCD
$$
where the maps related to the group $\Qm$ are indexed by $\Qm$.
\endproclaim

{\smc Definition 3.12.} The second four term exact sequence:
$$\eightpoint\CD
\thtw(H, \T)@>\Res_\Qm>> \La(\Hm, L, M, \T)@>\d_\Qm>>
\thasout(G, N, \T)@>\Inf_\Qm>> \ththr(H, \T)\\
\endCD
$$
will be called the {\it reduced} modified Huebschmann-Jones-Ratcliffe 
exact sequence
or simply the {\it reduced modified} HJR-exact sequence.

{\smc Remark 3.13.} The advantage of the reduced modified {\hjr} over 
the non-reduced one is that all the groups
involved are obviously compact, while the non-compactness of  the 
coefficient group $A$ in the non
reduced one forces us  to prove  the compactness of the cohomology 
group by examining the group of
cocycles and coboundaries.

\head{\bf\S4. Outer actions of a Countable Discrete Amenable  Group on an
AFD Factor of Type I\!I\!I$\pmb{_\la}$, $\pmb{0<\la<1}$.}
\endhead

We first apply the result of the last section to  the outer automorphism group
$\Out(\sM)$ by taking $\Out(\sM)$ as $G$ and $\thth1(\R,
\sU(L^\infty(\rt'z)))\cong \rtz$  as $N$:
\proclaim{Theorem 4.1} Suppose that $\sM$ is a separable factor of
type {\threel},
$0<\la<1$. The intrinsic modular invariant $\Obm(\sM)$ lives in the group{\rm:}
$$
\Obm(\sM)=([c], \nu_\sM)\in \thmsout(\Out(\sM), \thth1, \T),
$$
where $\thth1$ is the image $\{\dsigs: s\in \R\}$ of  the modular automorphism
group $\{\sigfs: s\in \R\},
\f\in \frak W_0(\sM),$ in the quotient group $\Out(\sM)=\Aut(\sM)/\Int(\sM)$,
$\nu_\sM$ is the identity map of $\thth1$ onto itself,
and $\txm$ is the modulus map $\txm\!:\da\in \Out(\sM)\mapsto
\mod(\da)\in\rt'z$.
The group
$\thmsout(\Out(\sM), \thth1, \T)$ is a non-separable compact abelian group.
\endproclaim
In view of \cite{KtT2: Theorem 3.2}, there is nothing left for the validity of
the assertion. But we want to identify the cocycle $c\in
\thsthr(\Outt(\tM)_\txm, \T)$
directly, where
$$
\Outt(\tM)_\txm=\{(p, s)\in \Outt(\tM)\times \R: \txm(p)=\ds\}.\tag4.1
$$

Before going further, let us fix notations for quotient maps:
$$\allowdisplaybreaks\aligned
&\hskip.5in \pi: \Out(\sM)\mapsto \Outt(\tM), \\
  \tpi&: \Aut(\sM)\mapsto \Outt(\tM)=\Aut(\sM)/\cntr(\sM),\\
\pi_0&: \Aut(\sM)\mapsto \Out(\sM)=\Aut(\sM)/\Int(\sM),\\
&\hskip.7in \tpi=\pi\scirc \pi_0.
\endaligned
$$
  Fix a generalized trace $\p\in \fW_0(\sM)$ so that $\sig_T^\p=\id$.
The  one-parameter unitary group $\{\p^{it}: t\in \R\}$ generates a
{\vna} $\sA^\p$
isomorphic  to $L^\infty(\R)$ and the non-commutative flow $\th$ restricted to
$\sA^\p$ is identified with the translation $\rho$:
$$
(\rho_tf)(x)=f(x+t), \quad f\in L^\infty(\R), x, t \in \R.
$$
We identify $\sA^\p$ and $L^\infty(\R)$ and $\p$ is then given by the function:
$$
\p(x)=e^{-x}, \quad x \in \R.
$$
The center $\sC$ of $\tM$ is then represented by the fixed point 
subalgebra $(\sA^\p)^{\th_{T'}}$ of
$\sA^\p$, i.e., the subalgebra of periodic functions with period $T'$.
We refer \cite{Tk2: Exercise XII.6, page 455 } for detail.

\proclaim{Lemma 4.2} If $\p\in \fW_0(\sM)$ is a generalized trace, 
i.e., a {\fwt} 
with period $T$
and $\p(1)=\infty$, then $\sM$ and $u_\p(\dst)=b_\p(\sz(\dst)), s\in 
\R,$ generates the discrete core
$\tM_d$.
\roster
\item"i)" The periodic one parameter unitary group $\{u_\p(\dst): 
s\in \R\}$ is represented by the
following function after $\sA^\p$ is identified with $\Linfr$\:
$$\aligned
u_\p(\ds; x)&=\exp\Big(\txti T'\bracettt{\ds}\Big[\frac
x{T'}\Big]\Big)=\Big\langle
\ds,  \Big[\frac x{T'}\Big]\Big\rangle,
\endaligned\quad x \in \R, \ds\in \rtz,\tag4.2
$$
which is also represented as a function of $\p$ in the following 
form\:
$$\aligned
u_\p(s)&=\exp\left(\txti T's\left[\frac {-\log\ \p}{T'}\right]\right)\\
&=\sum_{n\in\Z}\la^{-\txti ns}\chi_{(\la^{n+1}, \la^n]}(\p),\quad s\in
\R
\endaligned\tag4.2$'$
$$
where $\chi_{(\la^{n+1},\la^n]}(\p)$ 
means the spectral projection of $\p$ corresponding to 
the half open interval $(\la^{n+1},
\la^n]$.
\item"ii)" The cocycle $\sz(\ds)$ is cobounded by $u_\p(\ds)$ in $\sU(\sA^\p)$
relative to $\th${\rm:}
$$
\sz(\ds, t; x)=\Big\langle \ds,
\Big[\frac{x+t}{T'}\Big]-\Big[\frac{x}{T'}\Big]\Big\rangle, \quad \ds\in
\rtz,\  x \in\R;\tag4.3
$$
\item"iii)" $\Aut(\sM)_\txm$ acts on the discrete core $\tM_d$, i.e., 
if $(\a, s)\in \Aut(\sM)_\txm$,
then $\th_s\scirc\a$ leaves $\tM_d$ globally invariant.
\endroster
\endproclaim
\demo{Proof} The claims (i) and (ii) follow directly from \cite{FT2} 
with sign change in the
coboundary operation $\partth$. So we prove only (iii). First, 
observe that if $(\a, s)\in
\Aut(\sM)_\txm$, then $\th_s\scirc \a$ acts on the center $\sC$ 
trivially, i.e., it acts as the identity
since the actions $\mod(\a)$ and $\th_s$ cancel each other on $\sC$. 
Hence there exists $u\in
\sU(\sM)$ such that $\th_s\scirc \a(\p)=u\p u^*$. Thus we get
$$\allowdisplaybreaks\aligned
\th_s\scirc \a(u_\p(\dst))&=\th_s\scirc 
\a\left(\sum_{n\in\Z}\la^{-\txti ns}\chi_{(\la^{n+1},
\la^n]}(\p)\right)\\
&=\sum_{n\in\Z}\la^{-\txti ns}\chi_{(\la^{n+1},
\la^n]}(\th_s\scirc \a(\p))\quad\\\
&=\sum_{n\in\Z}\la^{-\txti ns}\chi_{(\la^{n+1}, \la^n]}(u\p u^*)\\
&=u\left(\sum_{n\in\Z}\la^{-\txti ns}\chi_{(\la^{n+1}, \la^n]}(\p)\right)u^*\\
&=uu_\p(\dst)u^*.
\endaligned
$$
Hence $\th_s\scirc \a(u_\p(\dst))\in \sM\vee \{u_\p(\dst): \dst\in 
\rtz\}\dprime=\tM_d$. This
completes the proof.
\qed
\enddemo
The generalized trace $\p$ gives rise to the semi-direct product decomposition:
$$
\cntr(\sM)=\Int(\sM)\rtimes_\sigp \left(\rtz\right).\tag4.4
$$
We are going to use the notation $\sig_\ds, \ds\in \rtz,$ for the 
element of $\cntr(\sM)$ corresponding
to an element $\ds=\dst=s+T\Z\in \rtz$.
\proclaim{Lemma 4.3} {\rm i)} It is possible to choose a 
cross-section $u$\: $m\in \cntr(\sM)\mapsto
u(m)\in
\tsU_0(\sM)\i \tM_d$ such that $u(m)\in \sU(\sM)$ if $m\in\Int(\sM)$ and
$$
u(\sig_\ds m)=u_\p(\ds)u(m), \quad \ds\in \rtz, \quad m\in \cntr(\sM).\tag4.5
$$

{\rm ii)} There exists a cross-section $g\in \Out(\sM)\mapsto \a_g\in 
\Aut(\sM)$ such that
\roster
\item "a)" Each $\a_g, g\in \Out(\sM)$, transforms $\p$ to a scalar 
multiple of $\p$, i.e., $\a_g(\p)$
and $\p$ are proportional. Consequently $\a_g$ and $\sigp$ commute;
\item"b)"
$$
\a_{\sig_\ds g}=\sig_\ds^\p\scirc \a_g, \quad \ds\in \rtz.\tag4.6
$$
\item"c)" The associated $\Int(\sM)$-valued cocycle $\etain$ has the property\:
$$
\etain(\sig_\ds g, \sig_\dt h)=\etain(g, h), \quad g, h\in \Out(\sM), 
\ds, \dt \in \rtz,\tag4.7
$$
so that it is the pullback $\pi^*(\etaq)$ of an $\Int(\sM)$-valued 
two cochain $\etaq \in
\tcatw(\Outt(\tM),
\Int(\sM))$.
\endroster
\endproclaim
\demo{Proof} i) First choose a Borel cross-section $u$: $m\in 
\Int(\sM)\mapsto u(m)\in\sU(\sM)$ of
the adjoint map $\Ad$: $v\in \sU(\sM)\mapsto \Ad(v)\in\Int(\sM)$. 
Then extend the cross-section by
setting:
$$
u(\sig_\ds m)=u_\p(\ds)u(m), \quad \ds\in \rtz, m\in\Int(\sM).
$$
This gives a cross-section with the desired property.

ii) For any $\a\in \Aut(\sM)$, $\a(\p)$ is another generalized trace 
on $\sM$. Hence there exists a
scalar $\la\in \R$ and a unitary $v\in \sU(\sM)$ such that 
$\Ad(v)\scirc\a(\p)=\la \p$. In fact, $\la$ is
can be chosen to be $e^{\bracett'{\mod(\a)}}$. Hence it is possible 
to a representative $\a_g$ of
$g\in \Out(\sM)$. With this in mind, we select first a cross-section 
$p\in \Outt(\tM)\mapsto \a_p\in
\Aut(\sM)$ of the quotient map $\tpi$: $\Aut(\sM)\mapsto 
\Outt(\tM)=\Aut(\sM)/\cntr(\sM)$ such that
\roster
\item"a)" The weights $\a_p(\p)$ and $\p$ are proportional, i.e., 
$\a_p(\p)=e^{\bracett'{\mod(p)}}\p$;
\item"b)" The quotient map $\pi_0$: $\Aut(\sM)\mapsto 
\Aut(\sM)/\Int(\sM)$ maps $\a_p$ exactly
into the cross-section image $\fs(p)\in \Out(\sM)$ of $p$ which has 
been fixed already.
\endroster
The cross-section $\a$ generates a $\cntr(\sM)$-valued two cocycle:
$$\allowdisplaybreaks\aligned
\etaa(p, q)&=\a_p\scirc\a_q\scirc \ainv{pq}\in \cntr(\sM), \quad p, q 
\in \Outt(\tM).
\endaligned\tag4.8
$$
Then we have
$$\aligned
\pi_0(\etaa(p, q))&=\pi_0(\a_p)\pi(\a_q)\pi_0 (\ainv{pq})=\fs(p)\fs(q)\fs(pq)\inv\\
&=\fns(p, q), \quad p, q \in \Outt(\tM).
\endaligned
$$
Therefore, the semi-direct product decomposition: 
$\cntr(\sM)=\Int(\sM)\rtimes_\sigp \rtz$ gives a
decomposition of $\etaa$:
$$
\etaa(p, q)=\fns(p, q)\eta_{\text {in}}(p, q), \quad p, q \in \Outt(\tM),
$$
where $ \fns(p, q)\in \thth1$ and  $\eta_{\text{in}}(p, q)\in \Int(\sM)$
 commutes with $\{\sigps;\ s\in \R\}$.
Now decomposing each $g\in
\Out(\sM)$ in the form:
$$\allowdisplaybreaks\aligned
g&=\fm_\thth1(g)\fs(\pi(g)),
\endaligned
$$
and writing $\ds(g)\in \rtz$ for $\fm_\thth1(g)\in \thth1$, we set
$$\aligned
\a_g&=\sig_{\ds(g)}^\p\scirc \a_{\pi(g)}\in \Aut_\p'(\sM)=\{\a\in 
\Aut(\sM): \a\scirc\sigps=\sigps
\scirc \a, s\in \R\}.\\
\endaligned
$$
Observing for each pair $g, h\in \Out(\sM)$ that
$$\allowdisplaybreaks\aligned
gh&=\fm_\thth1(g)\fs(\pi(g))\fm_\thth1(h)\fs(\pi(h))\\
&=\fm_\thth1(g)\fm_\thth1(h)\fs(\pi(g))\fs(\pi(h)),  \quad \text{as } 
\thth1\i \text{Center of
}\Out(\sM),\\
&=\fm_\thth1(g)\fm_\thth1(h)\fns(\pi(g), \pi(h))\fs(\pi(gh))\\
&=\fm_\thth1(gh)\fs(\pi(gh));\\
\fm_\thth1(g)&\fm_\thth1(h)=\fm_\thth1(gh)\fns(\pi(g), \pi(h))\inv,
\endaligned
$$
we compute:
$$\aligned
\a_g&\scirc \a_h=\sigpdsg\scirc \a_{\pi(g)}\scirc
\sigpdsh\scirc\a_{\pi(h)}=\sigpdsgph\scirc\a_{\pi(g)}\scirc \a_{\pi(h)}\\
&=\sigpdsgph\scirc \etaa(\pi(g), \pi(h))\scirc\a_{\pi(gh)}\\
&=\sig_{\ds(gh)}^\p\scirc\fns(\pi(g), \pi(h))\inv\scirc \etaa(\pi(g), 
\pi(h))\scirc \a_{\pi(gh)}\\
&=\eta_{\text{in}}(\pi(g), \pi(h))\scirc \sigpdsgh\scirc \a_{\pi(gh)}\\
&=\eta_{\text{in}}(\pi(g), \pi(h))\scirc \a_{gh}\\
\endaligned
$$
Therefore, the map $\a$: $g\in \Out(\sM)\mapsto \a_g\in 
\Aut_\p'(\sM)$ is indeed an outer action of
$\Out(\sM)$. Furthermore, $\etain(\pi(g), \pi(h))$ belongs to the 
group $\Int(\sM)\cap
\Autp'(\sM)$ which is given by the normalizer $N(\sM_\p)$ of the 
centralizer $\sM_\p$ of $\p$.
 From its construction, $\etain$ satisfies the requirement of the 
lemma. This completes the proof.
\qed
\enddemo

Before going into the last step, we need the following:
\proclaim{Lemma 4.4} {\rm i)}  If $\a\in \Autp'(\sM)$, then $\a$ 
leaves $\p$ relatively invariant, so that we
have
$\txk\in\Hom(\Autp'(\sM), \R)$ such that
$$
\a(\p)=e^{\txk(\a)}\p, \quad \a\in \Autp'(\sM),\tag4.9
$$

and $\a$ and $\th_{-\txk(\a)}$ agree on $\sA^\p$, in
particular on $u_\p(\ds), \ds\in\rtz${\rm;}

{\rm ii)} If $(\a, s)\in \Autp'(\sM)_\txm$, then
$$\allowdisplaybreaks\aligned
\Big((\a\scirc\th_s)(u_\p(\dt))u_\p(\dt)^*\Big)(x)&=\left\langle \dt,
\frac{s-\txk(\a)}{T'}\right\rangle\\
\endaligned\tag4.10
$$
in particular the left hand side is constant in $x\in \R$, i.e., the unitary
$u_\p(\dt),\allowmathbreak \dt\in \rtz,$ is an eigen operator of 
$\a\scirc \th_s$.

\endproclaim
\demo{Proof}
If $\a$ is in $\Autp'(\sM)$, then
its extension to $\tM$, still denoted by $\a$, leaves $\sA^\p$
globally invariant and
$$
(\a f)(x)=f(x-\txk(\a)), \quad f\in \sA^\p, x \in \R,
$$
as seen below:
$$
(\a \p)(x)=e^{\txk(\a)}\p(x)=e^{\txk(\a)}e^{-x}=e^{-(x-\txk(\a))}
=\p(x-\txk(\a)).
$$
The center $\sC$ of $\tM$ is generated by $\p^{\txti T}$, so that it
is identified with
the subalgebra of periodic functions with period $T'$, i.e.,
$\sC=(\sA^\p)^{\th_{T'}}$. By Lemma 4.2,
the periodic one parameter unitary group:
$$
\{u_\p(\dst)\!: \dst\in \rtz\}=\{b_\p(\sz(\dst))\!: \dst\in\rtz\}
$$
is represented in $\sA^\p$ by the following functions:
$$
u_\p(\dst; x)=\exp\Big(\txti T's\Big[\frac x{T'}\Big]\Big), \quad x, s \in \R,
$$
which together with $\sM$ generates the discrete core $\tM_d$. Since $\a\in
\Autp'(\sM)$ and $\th_{-\txk(\a)}$ both scales  the generator $\p$ of
$\sA^\p$ in the
same way, they agree on $\sA^\p$. Hence we get
$$\aligned
(\a (u_\p(\ds)))(x)=u_\p(\ds; x-\txk(\a))&\\
=\exp\Big(\txti T'\bracettt{\ds}\Big[\frac {x-\txk(\a)}{T'}\Big]\Big),&
\endaligned\qquad \ds\in \rtz.
$$
If $(\a, s)\in \Autp'(\sM)_\txm$, then
$$\allowdisplaybreaks\aligned
(\a\scirc\th_s)(u_\p(\dt))(x)&=u_\p(\dt; x+s-\txk(\a))\\
&=\exp\Big(\txti T'\bracettt{\dt}\Big[\frac {x+s-\txk(\a)}{T'}\Big]\Big)
\endaligned
$$
and
$$\allowdisplaybreaks\aligned
(\a\scirc\th_s)&(u_\p(\dt))u_\p(\dt)^*(x)=u_\p(\dt;
x+s-\txk(\a))\overline{u_\p(\dt;
x)}\\
&=\exp\Big(\txti T'\bracettt{\dt}\Big(\Big[\frac {x+s-\txk(\a)}{T'}\Big]
-\Big[\frac {x}{T'}\Big]\Big)\Big)\\
&=\exp\Big(\txti {T'}\bracettt{\dt}\frac{s-\txk(\a)}{T'}\Big)\\
&=\left\langle \dt, \dfrac{s-\txk(\a)}{T'}\right\rangle,
\endaligned
$$
where the last pairing makes a sense because $s-\txk(\a)\in {T'}\Z$ for $(\a, s)\in \Autp'(\sM)_\txm $.  This completes the proof.
\qed
\enddemo
Now the next lemma completes the proof of Theorem 4.1.
\proclaim{Lemma 4.5} If a cross-section $u$\: $m\in \cntr(\sM)\mapsto 
u(m)\in\tsU_0(\sM)$ is the one
given by {\rm Lemma 4.3}, then the natural choices of $u(g, h), g, 
h\in\Out(\sM)_\txm$, and  $u(\tp, \tq),
\tp=(p, s), \tq=(q, t) \in\Outt(\tM)_\txm$,
  by
$$\aligned
u(g, h)=u(\etain(g, h))\in \sU(\sM)\quad \text{and}\quad u(\tp, 
\tq)=u(\etaa(p, q))\in\tsU_0(\sM)
\endaligned
$$
give the following\:
\roster
\item"i)"
$$\allowdisplaybreaks\aligned
\a_g\scirc \a_h&=\Ad(u(g, h))\scirc \a_{gh};\\
\a_\tp\scirc \a_\tq&=\Ad(u(\tp, \tq))\scirc \a_{\tp\tq}.
\endaligned
$$
\item"ii)" The associated cocycles $c_\sM^{\text o}\in 
\tztw(\Out(\sM), \T)$ given by\:
$$
c_\sM^\circ(g, h, k)=\a_g(u(h, k))u(g, hk)\{u(g, h)u(gh, k)\}^*, \ g, 
h, k\in \Out(\sM),
$$
gives the cohomology class $[c_\sM^\circ]\in \ththr(\Out(\sM), \T)$ 
which is the intrinsic obstruction
$\Ob(\sM)$.
\item"iii)" The cocycle $c_\sM$ associated with the choice of $u$\:
$$
c_\sM(\tp, \tq, \tr)=\a_\tp(u(\tq, \tr))u(\tp, \tq\tr)\{u(\tp, 
\tq)u(\tp\tq, \tr)\}^*,\ \tp, \tq, \tr\in \Outt(\tM)_\txm
$$
is a standard coycle relative to the distinguished central element 
$z_0=(1, T')\in \Outt(\tM)_\txm$ and the
$d$-part is given by\:
$$
d(q, r)=\exp(\txti T'\bracett'{\fns(q, r)}), \quad q, r \in \Outt(\tM),
$$
so that
$$\allowdisplaybreaks\aligned
c_\sM(\tp z_0^m, \tq z_0^n, \tr z_0^\ell)=\big\langle \fns(q, r), 
m\big\rangle c_\sM(\tp, \tq, \tr),
\endaligned
$$
for each $m, n, \ell \in \Z$ and $\tp, \tq, \tr\in \Outt(\tM)_\txm$.
\item"iv)" The pair $(c_\sM, \id)$ belongs to 
$\tzsthr(\Outt(\sM)_\txm, \T)*_\fs \Hom_{\Out(M)}(\thth1, \thth1)$ and
its cohomology class $[c_\sM, \id]$ in $\tH_{\txs, 
\txm}^\out(\Out(\sM), \thth1, \T)$ is in fact the
intrinsic modular obstruction $\Obm(\sM)$.
\endroster
\endproclaim
\demo{Proof} We see immediately for each pair $g, h\in
\Out(\sM)$
$$\allowdisplaybreaks\aligned
\a_g\scirc \a_h&=\a_g\scirc \a_h=\etain(g, h)\scirc 
\a_{gh}=\Ad(u(\etain(g, h)))\scirc \a_{gh}\\
&=\Ad(u(g, h))\scirc\a_{gh},
\endaligned
$$
and for each pair $\tp=(p, s), \tq=(q, t)\in \Outt(\tM)_\txm$
$$\allowdisplaybreaks\aligned
\a_\tp\scirc \a_\tq&=\th_s\scirc \a_p\scirc \th_t\scirc 
\a_q=\th_{s+t}\scirc\a_p\scirc
\a_q=\th_{s+t}\scirc\etaa(p, q)\scirc \a_{pq}\\ &=\tAd(u(\etaa(p, 
q)))\scirc \a_{\tp\tq}=\tAd(u( p,
q))\scirc \a_{\tp\tq}.
\endaligned
$$
Consequently, $c_\sM^\circ$ and $c_\sM$ are both three cocycles and 
the former gives the intrinsic
obstruction $\Ob(\sM)$ in the cohomology group $\ththr(\Out(\sM), 
\T)$. To see the standardness of
$c_\sM$, we just examine:
$$\allowdisplaybreaks\align
c_\sM(\tp z_0^m&, \tq z_0^n, \tr z_0^\ell)=\a_\tp\scirc 
\th_{mT'}(u(\tq z_0^n, \tr z_0^\ell))u(\tp z_0^m, \tr
z_0^{n+\ell})\\
&\hskip1in\times
\{u(\tp z_0^m, \tq z_0^n)u(\tp\tq z_0^{m+n}, \tr z_0^\ell)\}^*\\
&=\a_\tp\scirc\th_{mT'}(u(q, r))u(p, qr)\{u(p, q)u(pq, r)\}^*\\
&=\a_\tp\scirc\th_{mT'}\Big(u_\p(\fns(q, r)) u(\etain(q, r))\Big)u(p, 
qr)\{u(p, q)u(pq, r)\}^*\\
&=\langle \fns(q, r), m\rangle\a_\tp\Big(u_\p(\fns(q, r)) u(\etain(q, 
r))\Big)\\
&\hskip1in \times
u(p, qr)\{u(p, q)u(pq, r)\}^*\\
&=\langle \fns(q, r), m\rangle\a_\tp\Big(u(q, r )\Big)
u(p, qr)\{u(p, q)u(pq, r)\}^*\\
&=\langle \fns(q, r), m\rangle c_\sM(\tp, \tq, \tr).
\endalign
$$

Since the $d$-part is given by the two cocycle $\fns$ itself, the 
pair $(c_\sM, \id)$ belongs to the fiber
product $\tzthr(\Outt(\sM)_\txm, \T)*_\fs\Hom_{\Out(M)}(\thth1, \thth1)$.
\qed
\enddemo

\proclaim{Theorem 4.6}{\rm i)} If G is a discrete group and $\a$ is 
an outer action of $G$
on a factor $\sM$ of type {\threel}, $0<\la<1$, then
the modulus  $\txm=\txm_\a\!: g\in G\mapsto \mod(g)\in \rtz'$ of $\a$, the
normal subgroup $N=N(\a)=\a\inv(\cntr(\sM))$, a homomorphism $\nu_\a$\:
$m\in N\mapsto
\da_m\in \thth1\cong\rtz$ and the ``pullback"
$[c_\a]=\a^*([c_\sM])\in\thmsout(G, N,
\T)$ of the intrinsic modular obstruction, to be termed the modular
obstruction of $\a$,
are outer conjugacy invariants of $\a$.

{\rm  ii)} If G is a countable discrete amenable group and the
factor $\sM$ of type
{\threel} is approximately finite dimensional, then the invariants
$\{\txm_\a, N(\a),
[c_\a], \nu_\a\}$ determines the outer conjugacy class of $\a$. The
group $\thmsout(G,
N, \T)$ is a separable compact abelian group.
\endproclaim

{\smc Remark 4.7.} The pullback in the theorem needs a 
qualification. As the cross-section $\fs:
\Outt(\tM)\mapsto \Out(\sM)$ is only guaranteed by the axiom of 
choice, we have no idea if it consistent
with the map $\da: G\mapsto \Out(\sM)$, for example it can happen 
that $\fs(\Outt(\tM))\cap
\da(G)=\{\id\}$. Namely, we cannot pull back the cross-section $\fs$ 
of $\Outt(\tM)$. So we have to work
with a cross-section $\fs: Q=G/N\mapsto G$ directly. But this does 
not change the picture concerning the
modular obstruction $\Obm(\a)$. If we consider all cross-sections of 
$\Outt(\tM)$ and form the group
$\tH^\out(\Out(\sM), \thth1, \T)$ as in \cite{KtT2: Page 218}, in 
which we locate the intrinsic modular
obstruction, then one can pull back $\Obm(\sM)$ to form the modular 
obstruction $\Obm(\a)\in
\tH^\out(G, N, \T)$ since cross-section of $Q\mapsto G$ can be 
carried to $\Outt(\tM)$ as a part of
a cross-section of $\Outt(\tM)$, i.e., we can have a cross-section 
$\fs$ of $\Outt(\tM)$ so that
$\da\inv(\fs(\Outt(\tM))N=G$, which enables us to pull back the 
structure concerning $\Aut(\sM),
\Out(\sM)$ and $\Outt(\tM)$.

\proclaim{Theorem 4.8} Suppose that $\sM$ is a factor of type {\threel},
$0<\la<1,$ and that $\p$ is a periodic {\fwt} with period $T=-2\pi/\log \la$.
Let $\txm\!:\a\in \Aut(\sM)\mapsto \mod(\a)\in\rt'z=\Aut_\th(\sC)$ be the
modulus homomorphism of $\Aut(\sM)$ to the
automorphism group $\Aut_\th(\sC)$ of the flow of weights $\{\sC, \R, \th\}
=\{L^\infty(\rt'z), \text{\rm Translation}\}$ which is identified with the
quotient group $\rt'z$, then the discrete core $\tM_d$ gives rise to an
$\Aut(\sM)_\txm$-equivariant commutative square of exact sequences{\rm:}
$$\allowdisplaybreaks\CD
@.1@.1@.1\\
@.@VVV@VVV@VVV\\
1@>>>\T@>>>\T@>>>1@>>>1\\@.@VVV@VVV@VVV\\
1@>>>\sU(\sM)@>>>\tsU_0(\sM)@>\part_{\th_{T'}}
>\underset{u_\p}\to\longleftarrow>\rtz@>>>0
\\@.@VVV@VVV@VVV\\
1@>>>\Int(\sM)@>>>\cntr(\sM)@>\dot\part_{\th_{T'}}>\underset{\sigp}\to\longleftarrow
>\rtz@>>>0\\ @.@VVV@VVV@VVV\\@.1@.1@.0
\endCD
$$
where $\tsU_0(\sM)$ is the unitary normalizer of $\sM$ in the 
discrete core $\tM_d$, i.e.,
$$
\tsU_0(\sM)=\tsU(\sM)\cap \tM_d=\{u\in\sU(\tM_d): u\sM u^*=\sM\}.
$$
\endproclaim
This is an easy consequence of the previous discussion, in particular 
Lemma 4.2.iii and
so we leave the proof to the reader.

{\smc Definition 4.9.} The above exact square is called the {\it 
reduced characteristic square} of $\sM$.

\head{\bf \S5. Outer actions of a Countable Discrete Amenable  Group on an
AFD Factor of Type I\!I\!I$\pmb{_1}$}
\endhead

The triviality of the flow of weights on a factor of type {\threeone} makes the
characteristic square very simple:
$$\CD
@.1@.1@.1\\
@.@VVV@VVV@VVV\\
1@>>>\T@>>>\T@>>>1@>>>1\\
@.@VVV@VVV@VVV\\
1@>>>\sU(\sM)@>>>\tsU(\sM)@>\partth>\underset{b_\f}\to\longleftarrow
>\R@>>>0\\
@.@V\Ad VV@V\tAd VV@|\\
1@>>>\Int(\sM)@>>>\cntr(\sM)@>\dot\part_\th>\underset{\sigf}\to\longleftarrow
>\R@>>>0\\ @.@VVV@VVV@VVV\\
@.1@.1@.0
\endCD
$$
Furthermore, the horizontal exact sequences split nicely. When we view
$\R$
as a central subgroup of $\Out(\sM)$, we denote it by
$\thth1$.
We will identify $\thth1$ and $\R$
frequently to avoid heavy notations in the case of type \threeone.

\proclaim{Theorem 5.1} Let $\sM$ be a factor of type {\threeone}. Fix a
cross-section
$$
\sp\!:p\in\Outt(\tM)\mapsto \sp(p)\in\Out(\sM)
$$
of the quotient map
$\pi\!: \Out(\sM)\mapsto \Outt(\tM)=\Out(\sM)/\thth1$ and the associated
$\R$-valued cocycle{\rm:}
$$
\fn_\R(p, q)=\fs(p)\fs(q)\fs(pq)\inv\in \R\cong\thth1,\quad p, q \in
\Outt(\tM).
$$
It is possible to select cross-sections $\a$\: $g\in \Out(\sM)\mapsto
\a_g\in \Aut(\sM)$  and $u\!: m\in \cntr(\sM)\mapsto u(m)\in\tsU(\sM)$ so that
\roster
\item"i)" $m=\tAd(u(m)), \quad m\in \cntr(\sM);$
\item"ii)" the associated intrinsic modular obstruction cocycle $c_\sM$ takes
the form{\rm:}
$$\aligned
c_\sM&(\tp, \tq, \tr)=\exp(-\txti s\fn_\R(q, r))c_{\Outt(\sM)}(p, q, r),
\endaligned
$$
for each $\tp=(p, s), \tq, \tr\in\Outt(\tM)\times \R$.
\endroster
Let $\a$ be an outer action of a {\cdg} $G$ on $\sM$ with
$N=\a\inv(\cntr(\sM))$ and $\nu_\a(m)=\da_m\in \R\cong\thth1, m \in N$. Fix
a cross-section $\fs$\: $Q=G/N\mapsto G$ of the quotient map $\pi: G\mapsto Q$
along with the associated $N$-valued cocycle{\rm:}
$$
\fnn(p, q)=\fs(p)\fs(q)\fs(pq)\inv\in N,\quad p, q \in Q.
$$
Then the pullback cocycle
$c_\a$ of $c_\sM$ by $\a^*$ takes the form{\rm:}
$$
c_\a(\tp, \tq, \tr)=\exp(-\txti s\nu_\a(\fnn(q, r)))c_Q(p, q, r),\quad
$$
for each $\tp=(p, s), \tq, \tr\in
\wtQ=Q\times \R$. Its cohomology class $([c_\a], \nu_\a)\in 
\thsout(G, N, \T)$ is the
modular obstruction $\Obm(\a)$ of $\a$.
\endproclaim

\demo{Proof} Fix a dominant weight $\f$ on $\sM$ and observe that
$$
\Aut(\sM)=\Int(\sM)\Autfm, \quad \Aut_\f(\sM)=\left\{\a\in \Aut(\sM): 
\f\scirc \a=\f\right\}.
$$
Then we have
$$\aligned
\Out(\sM)&=\Aut(\sM)/\Int(\sM)=\Autfm/(\Autfm\cap\Int(\sM)),\\
\Outt(\sM)&=\Aut(\sM)/\cntr(\sM)=\Autfm/(\Autfm\cap\cntr(\sM)),\\
&\qquad\tsU(\sM)=\sU(\sM)\rtimes_\sigf\R.
\endaligned
$$
The invariance $\f=\f\scirc \tAd(u), u\in \tsU(\sM)$, of $\f$ gives 
the decomposition $u=v\f^{\txti s}$ for some
$v\in \sU(\sM_\f)$ and $s\in \R$. Hence we get the decomposition:
$$
\tAd\inv(\cntr(\sM)\cap\Autfm)\cong\sU(\sM_\f)\times\R.
$$
Therefore, we can choose a cross-section $p\in\Outt(\sM)\mapsto 
\a_p\in \Autfm$ such that
$$
\a_p\scirc \a_q=\Ad(v(p, q))\scirc \sig_{\fn_\R(p, q)}^\f\scirc\a_{pq}, \quad v(p, 
q)\in \sU(\sM_\f), p, q \in \Outt(\sM).
$$
We then set
$$
\a_g=\sig_{\txs(g)}\scirc\a_{\pi(g)}, \quad g\in \Out(\sM),
$$
where
$$
g=\sig_{\txs(g)}^\f\sp(\pi(g)).
$$
Setting
$$
u(\tp, \tq)=v(p, q)\f^{\txti \fn_\R(p, q)}, \quad \tp=(p, s), \tq=(q, 
t) \in \Outt(\sM)\times \R,
$$
we compute the intrinsic modular obstruction cocycle:
$$\allowdisplaybreaks\aligned
c_\sM(\tp&, \tq, \tr)=\a_\tp(u(\tq, \tr))u(\tp\tq, \tr)\{u(\tp, 
\tq)u(\tp\tq, \tr)\}^*\\
&=\a_p\scirc\th_s(u(q, r))u(pq, r)\{u(p, q)u(pq, r)\}^*\\
&=\exp(-\txti s\fn_\R(q, r))\a_p(u(q, r))u(p, qr)\{u(p, q)u(pq, r)\}^*\\
&=\exp(-\txti s\fn_\R(q, r))\a_p(v(q, r))v(p, qr)\{v(p, q)v(pq, r)\}^*,
\endaligned
$$
where the last step follows from the fact that $\fn_\R$ is an 
$\R$-valued two cocycle over $\Outt(\sM)$. Therefore,
we conclude that
$$
c_\sM(\tp, \tq, \tr)=\exp(-\txti s\fn_\R(q, r))c_{\Outt(\sM)}(p, q, r).
$$
Now we look at $\a_g\scirc \a_h, g, h\in \Out(\sM)$:
$$\allowdisplaybreaks\aligned
\a_g\scirc \a_h&=\sig^\f_{\txs(g)}\scirc\a_{\sp(\pi(g))}\scirc 
\sig^\f_{\txs(h)}\scirc \a_{\sp(\pi(h))}\\
&=\sig^\f_{\txs(g)\txs(h)}\scirc \a_{\sp(\pi(h))}\scirc\a_{\sp(\pi(h))}\\
&=\sig^\f_{\txs(g)\txs(h)}\scirc\tAd(u(\pi(g), \pi(h)))\scirc 
\a_{\sp(\pi(gh))}\\
&=\sig^\f_{\txs(g)\txs(h)}\scirc\sig^\f_{\fn_\R(\pi(g), 
\pi(h))}\scirc\Ad(v(\pi(g), \pi(h))\scirc \a_{\sp(\pi(gh))}.
\endaligned
$$
Also we observe:
$$\allowdisplaybreaks\aligned
gh&=\txs(g)\sp(\pi(g))\txs(h)\sp(\pi(h))=\txs(g)\txs(h)\sp(\pi(g))\sp(\pi(h))\\
&=\txs(g)\txs(h)\fn_\thth1(\pi(g), \pi(h))\sp(\pi(gh))\\
&=\txs(gh)\sp(\pi(gh)), \quad g, h\in \Out(\sM);\\
&\txs(g)+\txs(h)=\txs(gh)-\fn_\R(\pi(g), \pi(h)).
\endaligned
$$
Plugging this into the previous computation, we get
$$
\a_g\scirc \a_h=\Ad(v(\pi(g), \pi(h))\scirc \a_{gh}, \quad g, h\in \Out(\sM).
$$
Therefore, the intrinsic obstruction cocycle $c$ is given by the pullback:
$$
c(g, h, k)=c_{\Outt(\sM)}(\pi(g), \pi(h), \pi(k)),\quad g, h, k\in\Out(\sM),
$$
of the restriction of the intrinsic modular obstruction cocycle 
$c_\tM$ to the subgroup $\Outt(\sM)$.
This completes the proof.
\qed
\enddemo

\bye